\documentclass[a4paper,12pt]{article}

\usepackage{latexsym}
\usepackage{amssymb}
\usepackage{theorem}
\usepackage{amsmath}
\usepackage{amscd}
\usepackage{graphicx}
\usepackage{xcolor}
\usepackage{url}

\newtheorem{theorem}{Theorem}[section]
\newtheorem{corollary}[theorem]{Corollary}
\newtheorem{lemma}[theorem]{Lemma}
\newtheorem{example}[theorem]{Example}
\newtheorem{proposition}[theorem]{Proposition}
\newtheorem{remark}[theorem]{Remark}
\newtheorem{definition}[theorem]{Definition}

\newcommand{\demo}{\par\noindent{\it Proof. \/}\ }
\newcommand{\enD}{\hfill $\Box$\vspace{3truemm} \par}
\newcommand{\R}{\mathbb{R}}

\newcommand{\bt}{\mbox{\boldmath $t$}}

\newcommand{\be}{\mbox{\boldmath $e$}}

\newcommand{\bv}{\mbox{\boldmath $v$}}
\newcommand{\bx}{\mbox{\boldmath $x$}}
\newcommand{\by}{\mbox{\boldmath $y$}}

\newcommand{\bgamma}{\mbox{\boldmath $\gamma$}}
\newcommand{\bnu}{\mbox{\boldmath $\nu$}}
\newcommand{\bmu}{\mbox{\boldmath $\mu$}}
\newcommand{\bE}{\mbox{\boldmath $E$}}

\begin{document}

\title{Horocyclic evolutes, parallels and  involutes of spacelike frontals in hyperbolic 2-space}

\author{Nozomi Nakatsuyama, Masatomo Takahashi and Anjie Zhou}

\date{\today}

\maketitle

\begin{abstract} 
The horocyclic evolutes of spacelike frontals in hyperbolic 2-space have already been defined. 
Using enveloid theorem, we now define the horocyclic parallel and involute  of a spacelike frontal  in hyperbolic 2-space as the normal envelopes of its normal and tangent horocycles, respectively. Meanwhile, we investigate the relations among horocyclic evolutes, parallels  and involutes.
\end{abstract}

\renewcommand{\thefootnote}{\fnsymbol{footnote}}

\section{Introduction}

In special relativity, the space of inertial frames is naturally hyperbolic, and connections between hyperbolic geometry and special relativity have long been investigated (cf. \cite{FDG2022,U2008}). 
In  \cite{FDG2022}, the geometrization of special relativity is achieved at the level of kinematic space by providing purely geometric definitions and descriptions for physical concepts and phenomena.
Among them, horocycles and horospheres are tangent to the ideal boundary of hyperbolic space.
Treating the horocycle or  horosphere as a ``straight line" or ``plane" yields    \textit{horospherical geometry}  (cf. \cite{BIR2011,I2009,IPRT2005}).
 Physically, horocycles and horospheres are objects where the energy and frequency of a photon remain constant. The relativistic Doppler effect arises when observers move across different horocycles or horospheres, providing  a   geometric basis for relativistic frequency shifts (cf. \cite{FDG2022}). Therefore, horocycles and horospheres are not only fundamental geometric objects in hyperbolic   space, but also directly explain real physical phenomena.
\par
On the other hand, the kinematic space of Newtonian mechanics is Euclidean space, while that of special relativity is hyperbolic space  (cf. \cite{FDG2022,V1910}).
Studying the geometric objects in these two spaces can  show the connection between classical mechanics and  special relativity.
 Parallels, evolutes and involutes  are all important objects that have been extensively studied   with significant applications (cf. \cite{EN2000,IPST2004,LT2023,NT2024,TG2023}). 
Evolutes and involutes can be defined by envelopes in Euclidean plane. The evolute of a regular plane curve is the envelope of its normal lines, while an involute is the normal envelope of its tangent lines.  Meanwhile, parallels of a regular plane curve can  be viewed as the normal envelopes of its normal lines. 
For curves that may contain singularities in   Euclidean plane, Fukunaga and  the second author  investigated  evolutes and  involutes of fronts and frontals in Euclidean plane  (cf. \cite{FT2013,FT2014,FT2015,FT2016}). 
In hyperbolic 2-space, 
authors in \cite{AII2015,CIT2024} have introduced the concepts of horocyclic evolutes of regular curves, spacelike fronts and   frontals.   However, the study of horocyclic involutes and  parallels   remains unexplored. In \cite{ZP}, the enveloid theorem in hyperbolic 2-space was provided by the third author and Pei. This motivates  the following definitions: the horocyclic involute of a spacelike frontal is the normal envelope of its  tangent horocycles, and the horocyclic parallel can be viewed as the normal envelope of  its normal horocycles.  
Based on this, we prove that the horocyclic evolutes, parallels   and involutes  are also spacelike frontals, investigating their properties and   relations.
The properties of evolutes can be discussed by distance squared functions and the
theories of Lagrangian, Legendrian singularity  (cf. \cite{A1990, A1995, AGV1986,BG1992,FT2014,IPST2004}).   A horocyclic evolute  or   involute  can  also be realized  as a discriminant set of a smooth function, respectively (Propositions \ref{prop3.5} and \ref{prop4.7}). 
\par
The paper is organized as follows. 
In Section \ref{S2}, we review preliminaries of spacelike frontals,  enveloid theorem  and criteria for singular points of curves in hyperbolic $2$-space. 
Section \ref{S3} is devoted to the study of horocyclic evolutes and   parallels. 
Since the horocyclic evolute and   parallel remain spacelike frontals, we  analyze their singularities and the relations between these singularities  using their curvatures (Proposition \ref{prop3.15}). 
We also show that the  horocyclic evolute of a spacelike frontal coincides with each   horocyclic evolute  of its horocyclic parallel  (Proposition  \ref{prop3.17}). 
In Section \ref{S4}, we define the horocyclic involute  of a spacelike frontal    using a Riccati equation. 
Moreover, we prove that the horocyclic evolute of a horocyclic involute of a spacelike frontal is the original curve, while the horocyclic involute of a horocyclic evolute is a horocyclic parallel of the original curve (Theorem \ref{th4.14}). In Section \ref{S5}, we present two examples.
\par
We shall assume throughout the whole paper that all maps and manifolds are $C^{\infty}$ unless the contrary is explicitly stated.

\bigskip
\noindent
{\bf Acknowledgement}. 
The first author is supported by JST SPRING (Grant Number JPMJSP2153).
The second author is partially supported by JSPS KAKENHI (Grant Number JP 24K06728).
The third author is  grateful for the financial support provided by the CSC-MuroranIT Scholarship (Grant Number 202506620020).

\section{Preliminaries}\label{S2}

Let $\R_1^3 $ denote \textit{Lorentz}-\textit{Minkowski} $3$-\textit{space} with the pseudo scalar product $\langle \bx,\by \rangle=-x_1y_1+x_2y_2+x_3y_3$ and the pseudo vector product 
$$
\bx \wedge \by =\det\begin{pmatrix}-\be_1&\be_2 &\be_3\\
x_1&x_2& x_3\\
y_1&y_2&y_3 
\end{pmatrix},
$$ 
where $\bx=(x_1,x_2,x_3)$, $\by=(y_1,y_2,y_3)$ and $\left\{\be_1,\be_2, \be_3\right\}$ is the canonical basis of $\R_1^3$. 
A vector $\bx \in\R_1^3 \backslash \left\{0\right\} $ is \textit{spacelike}, \textit{lightlike} or \textit{timelike} if $\langle \bx,\bx \rangle \textgreater 0$, $\langle \bx,\bx \rangle=0$ or $\langle \bx,\bx \rangle\textless0$, respectively. 
We define the \textit{hyperbolic} 2-\textit{space} as
$
H^2=\left\{ \bx\in\R_1^3 \mid \langle \bx,\bx\rangle=-1\right\}
$ and \textit{de Sitter} 2-\textit{space} as
$
S_1^2=\left\{ \bx\in\R_1^3 \mid \langle \bx,\bx\rangle=1\right\}.
$
It can be verified that
$\langle \bx,\by_1 \wedge \by_2  \rangle=\det(\bx,\by_1,\by_2),$ so $\by_1 \wedge \by_2$ and $\by_i$ are pseudo-orthogonal $ (i=1,2).$ 
We denote
$
\Delta_1=\left\{ (\bnu_1,\bnu_2) \in H^2\times S_1^2\mid \langle \bnu_1,\bnu_2 \rangle =0 \right\}. 
$ 
Given a lightlike constant vector $\ell$,  we define  a \textit{horocycle} by $HC(\ell,-1)=\left\{ \bx\in H^3 \mid \langle \bx,\ell\rangle=-1\right\}.$ 
In general, a horocycle is defined by $\langle \bx,\ell\rangle=c$, where  $c\neq0$.  However, if we choose $-\ell/c$  instead of $\ell$, then we have the above equation.

It is well known that there are two connected branches in hyperbolic 2-space. 
Here we only consider curves on the branch $\left\{ (x_1,x_2,x_3 ) \in {H} ^2 \mid  x_1>0\right\}$ and still denote it by $H^2$. 
In this paper, we denote $I,\widetilde{I},J,U$ as   intervals in $\R$.
\par
Let $\bgamma:I\rightarrow H^2$ be a regular curve, that is, $ \dot{\bgamma}(t)\neq0$ for all $t\in I$, where $\dot{\bgamma}(t)=(\mathrm{d}\bgamma/\mathrm{d}t)(t)$. 
Then we have the tangent vector $\bt(t)={\dot{\bgamma}(t)}/{|\dot{\bgamma}(t)|}$ and normal vector $\be(t)=\bgamma(t) \wedge \bt(t)$ of $\bgamma$. 
Therefore, $\{\bgamma(t),\bt(t),\be(t)\}$ is a moving frame along $\bgamma(t)$. 
We have the following hyperbolic Frenet-Serret formula: 
\begin{equation} \notag
\begin{pmatrix}
\dot{\bgamma}(t) \\
\dot{\bt}(t) \\
\dot{\be}(t) 
\end{pmatrix}
=
\begin{pmatrix}
0&|\dot{\bgamma}(t)|&0 \\
|\dot{\bgamma}(t)|&0&|\dot{\bgamma}(t)|\kappa_g(t) \\
0 &-|\dot{\bgamma}(t)|\kappa_g(t) &0
\end{pmatrix}
\begin{pmatrix}
\bgamma(t) \\
\bt(t) \\
\be(t) 
\end{pmatrix},
\end{equation}
where 
$
\kappa_g(t)={\det(\bgamma(t),\dot{\bgamma}(t),\ddot{\bgamma}(t))}/{|\dot{\bgamma}(t)|^3}
$ 
is the \textit{geodesic curvature} of $\bgamma$.
\begin{definition}[\cite{CT2016}] \rm
We call $(\bgamma,\bnu): I \rightarrow \Delta_1$ a \textit{spacelike Legendre curve} in $H^2$ if 	$\langle \dot{\bgamma}(t),\bnu(t) \rangle=0$ for all $t\in I$ and $\bgamma :I\rightarrow H^2$ is a \textit{spacelike frontal}.
Especially, $\bgamma $ is called a \textit{spacelike front} if $(\bgamma,\bnu )$ is an immersion.
\end{definition}
Let $\bmu(t) = \bgamma(t) \wedge \bnu(t) $. Then $\left\{\bgamma(t),\bnu(t),\bmu (t)\right\}$ is a moving frame along $\bgamma(t)$. 
We have the following Frenet type formula:
\begin{equation} \notag
\begin{pmatrix}
\dot{\bgamma}(t) \\
\dot{\bnu}(t) \\
\dot{\bmu}(t) 
\end{pmatrix}
=
\begin{pmatrix}
0&0&m(t) \\
0&0&n (t) \\
m(t) &-n(t) &0
\end{pmatrix}
\begin{pmatrix}
\bgamma(t) \\
\bnu(t) \\
\bmu(t)
\end{pmatrix}.
\end{equation}
The mapping $(m,n)$ is called the (\textit{spacelike Legendre}) \textit{curvature} of the spacelike Legendre curve $(\bgamma,\bnu)$. 
A point $t_0\in I$ is called an \textit{inflection point} of $\bgamma$ (or, $(\bgamma,\bnu)$) if $n(t_0)=0$.

\begin{definition}\rm 
Let $(\bgamma, \bnu), (\widetilde{\bgamma}, \widetilde{\bnu}) : I \to \Delta_1$ be spacelike Legendre curves. We say that $(\bgamma, \bnu)$ and $(\widetilde{\bgamma}, \widetilde{\bnu})$ are \emph{congruent as spacelike Legendre curves} if there exists  $A \in SO(1,2)$ such that
$
\widetilde{\bgamma}(t) = A(\bgamma(t))$,   $ \widetilde{\bnu}(t) = A(\bnu(t)) 
$
for all $t \in I$. 
Here, $$SO(1,2)=\left\{A\in M_3(\R)~\arrowvert~A^TGA = G, \;\det(A) = 1,\; G =
\begin{pmatrix}
-1 & 0 & 0  \\
0 & 1 & 0  \\
0 & 0 & 1 
\end{pmatrix}\right\} $$ and $A^T$  is the transpose matrix  of $A$.
\end{definition}

Then we have the following existence and uniqueness theorems in terms of the curvature of the spacelike Legendre curve.

\begin{theorem}[Existence theorem for spacelike Legendre curves  \cite{CPT2020}]
Let $(m, n) :$\\$ I \to \R^2 $ be a smooth mapping. Then there exists a spacelike Legendre curve $(\bgamma, \bnu) : I \to \Delta_1  $ whose associated curvature is $(m, n)$.
\end{theorem}

\begin{theorem}[Uniqueness  theorem for spacelike Legendre curves  \cite{CPT2020}]\label{th2.4}
 Let $(\bgamma, \bnu)$,\\ $(\widetilde{\bgamma}, \widetilde{\bnu}) : I \to \Delta_1$ be spacelike Legendre curves with curvatures $(m, n)$ and $(\widetilde{m}, \widetilde{n})$, respectively. Then $(\bgamma, \bnu)$ and $(\widetilde{\bgamma}, \widetilde{\bnu})$ are congruent as spacelike Legendre curves if and only if $(m, n)$ and $(\widetilde{m}, \widetilde{n})$ coincide.
\end{theorem}

\begin{lemma}[\cite{AII2015}]\label{lemma2.2}
Let $\left\{\bx_0,\bx_1,\bx_2\right\}$ be a   pseudo-orthonormal basis of $\R_1^3$, where $\bx_0\in H^2$, $\bx_1,\bx_2\in S_1^2$ and $\bx_2=\bx_0\wedge\bx_1$. 
The horocycle   tangent to $\bx_1$ at $\bx_0$ is given by 
$$
HC(\bx_0 \pm \bx_2,-1)=\left\{\bx\in H^2\mid\langle \bx,  \bx_0 \pm \bx_2 \rangle=-1\right\}
=\left\{\bx=\bx_0+s\bx_1+\frac{s^2}{2}(\bx_0\pm\bx_2) \mid s\in\R\right\}.
$$
\end{lemma}

\begin{definition}[\cite{CPT2020}] \rm
We call 
$(\bgamma,\bnu ):J\times I \rightarrow\Delta_1$ a \textit{one-parameter family of spacelike Legendre curves} in $H^2$ if $\langle\bgamma_{s},\bnu \rangle(s,t)=0$ for all $(s,t)\in J\times I$ and $\bgamma :J\times I \rightarrow H^2$ is a \textit{one-parameter family of spacelike frontals}, where ${\bgamma_s} (s,t)=(\partial\bgamma/\partial s)(s,t)$. 
\end{definition}
Let $\bmu(s,t) = \bgamma(s,t) \wedge \bnu(s,t)$. Then $\left\{\bgamma(s,t),\bnu(s,t),\bmu (s,t) \right\}$ is a moving frame along $\bgamma(s,t)$. 
Thus, we get
\begin{equation} \notag
\begin{aligned}
\begin{pmatrix}
\bgamma_{s}(s,t) \\
\bnu_{s} (s,t) \\
\bmu_{s}(s,t) 
\end{pmatrix}
&=
\begin{pmatrix}
0&0&m (s,t) \\
0&0&n (s,t) \\
m (s,t) &-n (s,t) &0
\end{pmatrix}
\begin{pmatrix}
\bgamma  (s,t)\\
\bnu  (s,t) \\
\bmu (s,t) 
\end{pmatrix},\\
\begin{pmatrix}
\bgamma_{t}(s,t) \\
\bnu_{t}  (s,t)\\
\bmu_{t} (s,t)
\end{pmatrix}
&=
\begin{pmatrix}
0&L (s,t) &M (s,t)\\
L (s,t) &0&N (s,t) \\
M (s,t) &-N (s,t) &0
\end{pmatrix}
\begin{pmatrix}
\bgamma (s,t) \\
\bnu (s,t) \\
\bmu (s,t)
\end{pmatrix}.
\end{aligned}
\end{equation}
The mapping $(m,n,L,M,N)$ is called the \textit{curvature} of $(\bgamma,\bnu)$.  More details are available from \cite{CPT2020}. 

\begin{definition}[\cite{ZP}]\label{def2.4}  \rm
Let  $(\bgamma,\bnu):J\times I \rightarrow\Delta_1$ be a one-parameter family of spacelike Legendre curves  with   curvature $(m,n,L,M,N)$ in $H^2$. 
Consider a curve $\be[\alpha]:U\rightarrow J\times I $,  $\be[\alpha](u)=(s(u), t(u))$, $u\in U$ and denote  {$\bE[\alpha](u)=\bgamma  \circ \be[\alpha](u)$},  where $\alpha\in[0,\pi)$. 
Then $\be[\alpha]$ is called a \textit{pre-$\alpha$-enveloid} and $\bE[\alpha]$ is an  \textit{$\alpha$-enveloid} of $\bgamma$ if   conditions (1) and (2)  are satisfied.
\begin{enumerate}
\item[$(1)$] The variability condition: $t$ is surjective and non-constant on any non-trivial subinterval of $U$.
\item[$(2)$] The $\alpha$-parallel condition: $\dot{\bE}[\alpha](u)$ and $\bmu \circ \be[\alpha](u)\cos\alpha+ \bnu \circ \be[\alpha](u)\sin\alpha$ are parallel for all $u\in U$.
\end{enumerate}
\end{definition}
\begin{lemma}[$\alpha$-enveloid theorem \cite{ZP}] \label{enveloidtheorem} 
Let $(\bgamma,\bnu ):J\times I \rightarrow\Delta_1$ be a one-parameter family of spacelike Legendre curves with curvature $(m,n,L,M,N)$. Suppose that
$\be[\alpha]:U\rightarrow J\times I $, $e[\alpha](u)=(s(u),t(u))$ is a smooth curve satisfying  variability condition in Definition \ref{def2.4}. Then $\be[\alpha]$ is a pre-$\alpha$-enveloid of $\bgamma $ if and only if
$$ 
\dot{s}(u)m \circ \be[\alpha](u)\sin\alpha+ \dot{t}(u)M \circ \be[\alpha](u)\sin\alpha- \dot{t}(u)L \circ \be[\alpha](u)\cos\alpha=0 
$$ for all $u\in U$.
\end{lemma}

We say that a smooth curve $\bgamma:(I,t_0)\rightarrow H^2$ at $t_0$ has a $(i,j)$-cusp, where $ (i, j)= (2, 3),(2, 5)$, $(3, 4),(3, 5) $  if $\bgamma$ is $ \mathcal{A} $-equivalent to the germ  $t \mapsto(t^i,t^j)$ at the origin. For curves with $(i,j)$-cusps in $\R^2$, the  criteria
are available  in \cite{BG1982,HS2019,P2001}. In order to consider singular points, we provide criteria of singular points of curves in $H^2$.
\begin{lemma}\label{lemma2.5}
Let $\bgamma:I\rightarrow H^2$ be a smooth curve and $t_0\in I$. Then the following assertions hold.
\begin{enumerate}
\item[$(1)$]  $\bgamma$ has a $(2,3)$-cusp at $t_0$ if and only if $\dot{\bgamma}(t_0)=0$ and $\det(\bgamma,\ddot{\bgamma},\dddot{\bgamma})(t_0)\neq0$.
\par
\item[$(2)$]  $\bgamma$ has a $(3,4)$-cusp at $t_0$ if and only if $\dot{\bgamma}(t_0)=0$, $\ddot{\bgamma}(t_0)=0$ and $\det(\bgamma,\dddot{\bgamma},\bgamma^{(4)})(t_0)\neq0$.
\par
\item[$(3)$] $\bgamma$ has a $(2,5)$-cusp at $t_0$ if and only if $\dot{\bgamma}(t_0)=0$, $\ddot{\bgamma}(t_0)\neq0$, $\dddot{\bgamma}(t_0)=c\ddot{\bgamma}(t_0)$ for some constant $c\in\R$  and $\det(\bgamma,\ddot{\bgamma},3\bgamma^{(5)}-10c\bgamma^{(4)})(t_0)\neq0$.
\par
\item[$(4)$] $\bgamma$ has a $(3,5)$-cusp at $t_0$ if and only if $\dot{\bgamma}(t_0)=0$, $\ddot{\bgamma}(t_0)=0$, $\det(\bgamma,\dddot{\bgamma},\bgamma^{(4)})(t_0)=0$ and  $\det(\bgamma,\dddot{\bgamma},\bgamma^{(5)})(t_0)\neq0$.
\end{enumerate} 
\end{lemma}
\demo 
Consider the diffeomorphism ${\pi}:  H^2\rightarrow D^2$ defined by
${\pi}(x,y,z)=\left({y}/({x +1}),{z}/({x +1}) \right),$
where $D^2=\left\{(a,b)\in\R^2\mid a^2+b^2<1\right\}$ is the Poincar\'{e} 2-disc. Let $\bgamma:I\rightarrow H^2$ be given by  $\bgamma(t)=(x(t),y(t),z(t))$ and $\overline{\bgamma}(t)=\pi\circ\bgamma(t)=\left({y(t)}/({x(t) +1}),{z(t)}/({x (t)+1}) \right),$ where $x,y,z:I\rightarrow \R$ are smooth functions.
Since $-x^2(t)+y^2(t)+z^2(t)=-1$ and $-x(t)\dot{x}(t)+y(t)\dot{y} (t)+z(t)\dot{z}(t)=0$ hold for all $t\in I$, we obtain that $\dot{\bgamma}(t )=(\dot{x}(t ),\dot{y}(t ),\dot{z}(t))=(0,0,0)$ if and only if $$\dot{\overline{\bgamma}}(t )=\left(\frac{\dot{y}(t )(x(t)+1)-y(t)\dot{x}(t)}{(x(t) +1)^2}, \frac{\dot{z}(t)(x(t)+1)-z(t)\dot{x}(t)}{(x (t)+1)^2} \right) =(0,0). $$ If $t_0$ is a singluar point of $\bgamma$, that is, $\dot{\bgamma}(t_0)=0$, then
  $\det(\bgamma,\ddot{\bgamma},\dddot{\bgamma})(t_0)=(x(t_0)+1)^2\det(\ddot{\overline{\bgamma}},\dddot{\overline{\bgamma}})(t_0)$. Therefore, $\bgamma$ has a $(2,3)$-cusp at $t_0$ if and only if $\det(\bgamma,\ddot{\bgamma},\dddot{\bgamma})(t_0)\neq0$ by the criteria of $(2,3)$-cusp in \cite{P2001}.

Furthermore, one can show that $\ddot{\bgamma}(t_0)=0$ if and only if $\ddot{\overline{\bgamma}}(t_0)=0$, and that $\dddot{\bgamma}(t_0)=c\ddot{\bgamma}(t_0)$ if and only if $\dddot{\overline{\bgamma}}(t_0)=c\ddot{\overline{\bgamma}}(t_0)$ for some constant $c\in\R$. 
If $\dot{\bgamma}(t_0)=\ddot{\bgamma}(t_0)=0$, we have 
$$\begin{aligned}
\det(\bgamma,\dddot{\bgamma},\bgamma^{(4)})(t_0)&=(x(t_0)+1)^2\det(\dddot{\overline{\bgamma}}, {\overline{\bgamma}}^{(4)})(t_0),\\ \det(\bgamma, {\bgamma}^{(4)},\bgamma^{(5)})(t_0)&=(x(t_0)+1)^2\det( {\overline{\bgamma}}^{(4)}, {\overline{\bgamma}}^{(5)})(t_0).\end{aligned}$$  It follows that the assertions $(2)$ and $(4)$ hold.
If $\dot{\bgamma}(t_0)=0$, $\ddot{\bgamma}(t_0)\neq0$ and  $\dddot{\bgamma}(t_0)=c\ddot{\bgamma}(t_0)$ for some constant $c\in\R$,  then
$\det(\bgamma,\ddot{\bgamma},3\bgamma^{(5)}-10c\bgamma^{(4)})(t_0) =(x(t_0)+1)^2\det( \ddot{\overline{\bgamma}},3\overline{\bgamma}^{(5)}-10c\overline{\bgamma}^{(4)})(t_0).$ 
This yields assertion $(3)$.
\enD

\begin{corollary}\label{cor2.6}
Let $(\bgamma,\bnu):I\rightarrow \Delta_1$ be a spacelike Legendre curve with curvature $(m,n)$.  Then we have the following.
\begin{enumerate}
\item[$(1)$]  $\bgamma$ is singular at $t_0$  if and only if $m(t_0)=0$.
\item[$(2)$]  $\bgamma$ has a $(2,3)$-cusp at $t_0$ if and only if $m(t_0)=0$, $ n(t_0)\neq0$, $\dot{m}(t_0)\neq0$.
\item[$(3)$] $\bgamma$ has a $(3,4)$-cusp at $t_0$ if and only if  $m(t_0)=\dot{m}(t_0)=0$, $ n(t_0)\neq0$, $\ddot{m}(t_0)\neq0$.
\item[$(4)$] $\bgamma$ has a $(2,5)$-cusp at $t_0$ if and only if  $m(t_0)=n(t_0)=0$, $\dot{m}(t_0)\neq0$, $\dot{m}(t_0)\ddot{n}(t_0)-\ddot{m}(t_0)\dot{n}(t_0)\neq0$.
\item[$(5)$] $\bgamma$ has a $(3,5)$-cusp at $t_0$ if and only if  $m(t_0)=n(t_0)=\dot{m}(t_0)=0$, $\ddot{m}(t_0)\neq0$, $\dot{n}(t_0)\neq0$.
\end{enumerate} 
\end{corollary}

\section{Horocyclic evolutes and   parallels of spacelike frontals in $H^2$}\label{S3}

Let $(\bgamma,\bnu ):I\rightarrow\Delta_1$ be a spacelike Legendre curve  with curvature $(m,n)$. Assume that there exists a smooth function $f:I\rightarrow\R  $ such that $m(t)+f(t)n(t)=0$ for all $t\in I$. 
Then the \textit{horocyclic evolute}  of $\bgamma$ (or, \textit{horocyclic evolute}  of $(\bgamma,\bnu)$)  is defined as (\cite{CIT2024}) 
\begin{equation}\notag
Ev^\pm(\bgamma)(t)=\bgamma(t)+f(t)\bnu(t)+\frac{f^2(t)}{2}(\bgamma(t)\pm\bmu(t)).
\end{equation} 

\begin{proposition}[\cite{CIT2024,FT2015}]
We denote  $\mathrm{Reg}(\bnu)=\left\{t\in I\mid n(t)\neq0\right\}$ as the set of regular points of  $\bnu$. Suppose that there exists a smooth function $f:I\rightarrow\R  $ such that $f(t)=-m(t)/n(t)$ on $\mathrm{Reg}(\bnu)$. Then   $f$ is
 unique if and only if $\mathrm{Reg}(\bnu) $ is   dense  in $ I $. 
\end{proposition}

The uniqueness condition of $f$ is a topological condition. In
this paper, we assume that $\mathrm{Reg}(\bnu)=\left\{t\in I\mid n(t)\neq0\right\}$  is   dense  in $ I $. 

\begin{remark}[\cite{AII2015}]\rm
The \textit{horocyclic evolute} of a regular curve $\bgamma:I\rightarrow H^2$ with nowhere zero geodesic curvature $\kappa_g$ is defined as 
$$\mathcal{E}v^\pm(\bgamma)(t)=\bgamma(t)+\frac{1}{\kappa_g(t)}\be(t)+\frac{1}{2\kappa^2_g(t)}(\bgamma(t)\pm\bt(t)).$$ 
\end{remark}

Let $(\bgamma,\bnu):I\rightarrow \Delta_1$ be a spacelike Legendre curve with curvature $(m,n)$ and
take the frame
$ 
\left\{\bx_0(t),\bx_1(t), \bx_2(t)\right\}=\left\{\bgamma(t),\bnu (t),\bmu(t) \right\} 
$. 
By Lemma \ref{lemma2.2}, we obtain the one-parameter family of normal horocycles of $\bgamma$  
\begin{equation}\notag
\bx_N^\pm(s,t)=\bgamma(t)+s\bnu(t)+\frac{s^2}{2}(\bgamma(t)\pm\bmu(t)). 
\end{equation}

It can be proved that $\bx_N^\pm $ is a one-parameter family of spacelike frontals.
\begin{proposition}\label{p0}
Under the above notations, $(\bx_N^\pm,\bnu_N^\pm):\R\times I\rightarrow \Delta_1$ is a one-parameter family of spacelike Legendre curves with curvature 
$$
\begin{aligned}
&\left({m}_N^\pm(s,t),{n}_N^\pm(s,t),{L}_N^\pm(s,t),{M}_N^\pm(s,t),{N}_N^\pm(s,t)\right)\\=&\left(\pm1,\pm1, \mp(m(t)+sn(t)),s^2n(t)/2+sm(t),s^2n(t)/2+sm(t)+n(t)\right),	
\end{aligned}
$$
where 
$
{\bnu}_N^\pm(s,t)=\mp\bmu(t)+s\bnu(t)+(s^2/2)(\bgamma(t)\pm\bmu(t)).
$
\end{proposition}
\demo  
By a direct calculation, 
$\bx_{Ns}^\pm(s,t)=\bnu(t)+s(\bgamma(t)\pm\bmu(t))$. Take $
{\bnu}_N^\pm(s,t)=\mp\bmu(t)+s\bnu(t)+(s^2/2)(\bgamma(t)\pm\bmu(t)),
$  we have  $ 
\langle  \bx_N^\pm(s,t),\bnu_N^\pm(s,t)\rangle=\langle   \bx_{Ns}^\pm(s,t),\bnu_N^\pm(s,t)\rangle=0
$ 
for all $(s,t)\in \R\times I$.   Hence $( \bx_N^\pm,\bnu_N^\pm)$ is a    one-parameter family of spacelike Legendre curves. Define 
$ \bmu_N^\pm(s,t)=\bx_N^\pm(s,t)\wedge\bnu_N^\pm(s,t)=s\bmu(t)\pm\bnu(t)\pm s\bgamma(t).$  Then $\{\bx_N^\pm(s,t),\bnu_N^\pm(s,t),\bmu_N^\pm(s,t)\}$ is a   moving frame along $\bx_N^\pm(s,t)$. Since   $$\begin{aligned}
	\bx_{Ns}^\pm(s,t)&=\bnu_{Ns}^\pm(s,t)=\bnu(t)+s(\bgamma(t)\pm\bmu(t))=\pm\bmu_N^\pm(s,t),\\
	\bx_{Nt}^\pm(s,t)&=\pm(s^2/2)m(t)\bgamma(t)\mp(s^2/2)n(t)\bnu(t)+((1+s^2/2)m(t)+sn(t))\bmu(t)\\&=\mp(m(t)+sn(t))\bnu_N^\pm(s,t)+(s^2n(t)/2+sm(t))\bmu_N^\pm(s,t),\\
\bv_{Nt}^\pm(s,t)&=\mp(1-s^2/2)m(t)\bgamma(t)\pm(1-s^2/2)n(t)\bnu(t)+(s^2m(t)/2+sn(t))\bmu(t)\\&=\mp(m(t)+sn(t))\bx_N^\pm(s,t)+(s^2n(t)/2+sm(t)+n(t))\bmu_N^\pm(s,t),\end{aligned}$$
    the curvature   of one-parameter family of spacelike Legendre curves $( \bx_N^\pm,\bnu_N^\pm)$ is given by \\ $\left(\pm1,\pm1, \mp(m(t)+sn(t)),s^2n(t)/2+sm(t),s^2n(t)/2+sm(t)+n(t)\right).$ 
\enD
Next, we show that the horocyclic evolute of $\bgamma$ is an envelope of the one‑parameter family of normal horocycles of $\bgamma$.
\begin{proposition}
Let $(\bgamma,\bnu ):I\rightarrow\Delta_1$ be a spacelike Legendre curve  with curvature $(m,n)$ 
and $(\bx_N^\pm,{\bnu}_N^\pm):\R\times I\rightarrow \Delta_1$ is a one-parameter family of spacelike Legendre curves defined by Proposition \ref{p0}. Suppose  that there exists a smooth function $f:I\rightarrow\R  $ such that $m(t)+f(t)n(t)=0$ for all $t\in I$.  
Then the  horocyclic evolute   of $\bgamma$,
$$
Ev^\pm(\bgamma)(t)=\bx_N^\pm\circ\be_N[0](t)=\bE_N^\pm[0](t) 
$$ is an envelope of $\bx_N^\pm$, where $\be_N[0]:I\rightarrow \R\times I$,
$ \be_N[0](t)=(f(t),t)$ is a pre-envelope of $\bx_N^\pm$.
\end{proposition}
\demo
Since $m(t)+f(t)n(t)=0$ holds for all $t\in I$,  we can obtain  $L_N^\pm(f(t),t)=\mp(m(t)+f(t)n(t))=0.$  Hence the curvature  of $(\bx_N^\pm,\bnu_N^\pm)$ satisfies
$$\dot{f}(t)m_N^\pm(f(t),t)\sin 0+M_N^\pm(f(t),t)\sin 0-L_N^\pm(f(t),t)\cos 0=0 $$ for all $t\in I$.
It follows from 
Lemma \ref{enveloidtheorem} that $\be_N[0] (t)=(f(t),t)$ is a pre-envelope of $ \bx_N^\pm $. 
Moreover,
\begin{equation}\notag
\begin{aligned}
Ev^\pm(\bgamma)(t)&= \bgamma(t)+f(t)\bnu(t)+\frac{f^2(t)}{2}(\bgamma(t)\pm\bmu(t))= {\bx}_N^\pm(f(t),t)= \bx_N^\pm\circ\be_N[0] (t)  =\bE_N^\pm[0](t)
\end{aligned}
\end{equation}
is an envelope of $\bx_N^\pm$.\enD

Let $F:\left(\R\times\R^r, (t_0,\bx_0)\right)\rightarrow\R$  be a smooth function germ. The \textit{discriminant set} of $F$ is  
$$
\mathcal{D}_F=\left\{ \bx \in \R^r \mid \ \text{there exists}\ t \text{ with } F(t,\bx)=\frac{\partial F}{\partial t}(t,\bx)=0 \right\}.
$$
Define a smooth  function $H_E^\pm:I\times H^2\rightarrow\R$, $H_E^\pm(t,\bx)=\langle\bgamma(t)\pm\bmu(t),\bx\rangle+1$. We have the following conclusion.
\begin{proposition}\label{prop3.5}
Under the above notations, suppose that there exists a smooth function $f:I\rightarrow\R  $ such that $m(t)+f(t)n(t)=0$ for all $t\in I$. 
Then the image of horocyclic evolute of $\bgamma$ is the discriminant set of $ H_E^\pm$. 
\end{proposition}
\demo For any $(t, \bx)\in I\times H^2$, there exist $a,b,c\in\R$ such that $\bx=a\bgamma(t)+b
\bmu(t)+c\bnu(t)$.
Under the condition $H_E^\pm(t,\bx)=\langle\bgamma(t)\pm\bmu(t),\bx\rangle+1=0$, we have $a=1\pm b$. Since $\langle\bx,\bx\rangle=-a^2+b^2+c^2=-(1\pm b)^2+b^2+c^2=-1$, it follows that $b=\pm c^2/2$. Thus, $$\bx=\left(1+\frac{c^2}{2}\right)\bgamma(t)+c\bnu(t)\pm\frac{c^2}{2}\bmu(t).$$
Differentiating $ H_E^\pm(t,\bx)$ with respect to 
$t$ and substituting $\bx$ yields 
$$ \frac{\partial H_{E}^\pm}{\partial t}(t,\bx)=\left\langle m(t)\bmu(t)\pm(m(t)\bgamma(t)-n(t)\bnu(t)),\bx\right\rangle=\mp(m(t)+cn(t)).
$$
By $(\partial H_E^\pm/\partial t)(t,\bx)=0$, we deduce $m(t)+cn(t)=0$.  Given the existence of a smooth function $f:I\rightarrow\R$ satisfying $m(t)+f(t)n(t)=0$ for all $t\in I$ and the assumption that $\mathrm{Reg}(\bnu)=\left\{t\in I\mid n(t)\neq0\right\}$  is   dense  in $ I $, we   conclude that  $c=f(t)$ holds. This completes the proof.
\enD

\begin{proposition}[\cite{CIT2024}]\label{prop3.1}
Suppose that there exists a smooth function $f:I\rightarrow\R  $ such that $m(t)+f(t)n(t)=0$ for all $t\in I$.  Then  
$\left({Ev}^\pm(\bgamma ), \bnu^\pm_E \right):I\rightarrow \Delta_1$ is a spacelike Legendre curve with curvature 
$$
\left(m^\pm_E(t), n^\pm_E(t)\right)= \left( - \dot{f}(t)\pm\frac{f^2(t)n(t)}{2},\pm \dot{f}(t)-\frac{f^2(t)n(t)}{2}+n(t)\right),
$$ 
where 
$$
\bnu^\pm_E(t)=\mp\frac{f^2(t)}{2}\bgamma(t)\mp f(t)\bnu(t)+ \left(1-\frac{f^2(t)}{2}\right)\bmu(t).
$$
\end{proposition}

Combining Corollary \ref{cor2.6} and Proposition \ref{prop3.1}, we have the following.
\begin{corollary}
Let $(\bgamma,\bnu ):I\rightarrow\Delta_1$ be a spacelike Legendre curve  with curvature $(m,n)$.	   Suppose that there exists a smooth function $f:I\rightarrow\R  $ such that $m(t)+f(t)n(t)=0$ for all $t\in I$.  For any $t_0\in I$, the following assertions hold.
\begin{enumerate}
\item[$(1)$]  $Ev^\pm(\bgamma) $ is singular at $t_0$ if and only if  $ - 2\dot{f}(t_0)\pm f^2(t_0)n(t_0)=0 $.

\item[$(2)$]  $Ev^\pm(\bgamma) $ has a $(2,3)$-cusp at $t_0$ if and only if  $n(t_0)\neq0$, $2(\dot{m}n-m\dot{n})(t_0)\pm m^2(t_0)n(t_0)=0$ and $2(\ddot{m}n-m\ddot{n})(t_0)\pm3m^2(t_0)\dot{n}(t_0)-m^3(t_0)n(t_0)\neq0$.

\item[$(3)$] $Ev^\pm(\bgamma) $ has a  $(3,4)$-cusp at $t_0$ if and only if  $n(t_0)\neq0$, $2(\dot{m}n-m\dot{n})(t_0)\pm m^2(t_0)n(t_0)=0$,   $2(\ddot{m}n-m\ddot{n})(t_0)\pm3m^2(t_0)\dot{n}(t_0)-m^3(t_0)n(t_0)=0$ and $ \mp m^2(t_0)\left(6\dot{n}^2+ 4n\ddot{n}-3m^2n^2 \right)(t_0)+4n(t_0)(m\dddot{n}-\dddot{m}n )(t_0)\neq0$.

\item[$(4)$] $Ev^\pm(\bgamma) $ has a  $(2,5)$-cusp at $t_0$ if and only if  $n(t_0)=\dot{f}(t_0)=0$, $ -2\ddot{f}(t_0)\pm f^2(t_0)\dot{n}(t_0)\neq0 $ and $ -\ddot{f}(t_0)\ddot{n}(t_0)+\dddot{f}(t_0)\dot{n}(t_0)\neq0 $.

\item[$(5)$] $Ev^\pm(\bgamma) $ has a  $(3,5)$-cusp at $t_0$ if and only if   $n(t_0)=\dot{f}(t_0)=  -2\ddot{f}(t_0)\pm f^2(t_0)\dot{n}(t_0)=0 $, $\dot{n}(t_0)\neq0$ and $-2\dddot{f}(t_0)\pm f^2(t_0)\ddot{n}(t_0)\neq0$. 
\end{enumerate}
\end{corollary}

\begin{corollary}
Let $(\bgamma,\bnu ):I\rightarrow\Delta_1$ be a spacelike Legendre curve  with curvature $(m,n)$.	   Suppose that there exists a smooth function $f:I\rightarrow\R  $ such that $m(t)+f(t)n(t)=0$ for all $t\in I$.  For any $t_0\in I$, the following assertions hold.
\begin{enumerate}
\item[$(1)$]   If $t_0$ is a $(2,3)$-cusp, $(2,5)$-cusp or $(3,5)$-cusp of $\bgamma$, then $t_0$ is a regular point of $Ev^\pm(\bgamma)$. 

\item[$(2)$] If  $t_0$ is a $(3,4)$-cusp of $\bgamma$, then  $t_0$ is a $(2,3)$-cusp of $Ev^\pm(\bgamma)$. 

\item[$(3)$]   If $t_0$ is a   $(2,5)$-cusp or $(3,5)$-cusp of $Ev^\pm(\bgamma)$, then $t_0$ is a singular point of $\bgamma$.
\end{enumerate}
\end{corollary}

It can be proved that $(\bgamma,-\bnu):I\rightarrow \Delta_1$ is a spacelike Legendre curve with curvature $(-m,n)$. 
Suppose that there exists a smooth function $f:I\rightarrow\R  $ such that $m(t)+f(t)n(t)=0$ for all $t\in I$. 
The existence of a horocyclic evolute  of $(\bgamma,\bnu)$ implies that of a horocyclic evolute  of $(\bgamma,-\bnu)$.
We denote the horocyclic evolute of $(\bgamma,-\bnu)$ by $\overline{Ev}^\pm(\bgamma)$ in this paper and have the following conclusion.  
\begin{proposition}\label{prop3.7}
$(\overline{Ev}^\pm(\bgamma),\overline{\bnu}_E^\pm):I\rightarrow\Delta_1$ is a spacelike Legendre curve with curvature 
$$
\left(\overline{m}_E^\pm(t),\overline{n}_E^\pm(t)\right)=\left(-\dot{f}(t)\mp\frac{f^2(t)n(t)}{2},\mp \dot{f}(t)-\frac{f^2(t)n(t)}{2}+n(t)\right),
$$
where $\overline{Ev}^\pm(\bgamma)=Ev^\mp(\bgamma)$
and $\overline{\bnu}_E^\pm=\bnu_E^\mp$.
\end{proposition}
\demo Since $\bgamma(t)\wedge(-\bnu(t))=-\bmu(t)$ and $(-f(t))n(t)-m(t)=-\left(f(t)n(t)+m(t)\right)=0$ for all $t\in I$,
the horocyclic evolute of $(\bgamma,-\bnu)$ is given by
$$\begin{aligned}
\overline{Ev}^\pm(\bgamma)(t)&=\bgamma(t)+(-f(t))(-\bnu(t))+\frac{(-f(t))^2}{2}(\bgamma(t)\mp\bmu(t))\\
&=\bgamma(t)+f(t)\bnu(t)+\frac{f^2(t)}{2}(\bgamma(t)\mp\bmu(t)) =Ev^\mp(\bgamma)(t).
\end{aligned}$$
Take $\overline{\bnu}_E^\pm =\bnu_E^\mp $. By Proposition \ref{prop3.1}, $(\overline{Ev}^\pm(\bgamma),\overline{\bnu}_E^\pm)=(Ev^\mp(\bgamma),\bnu_E^\mp)$ is a spacelike Legendre curve with curvature  $$\left(\overline{m}_E^\pm(t),\overline{n}_E^\pm(t)\right)=\left(m_E^\mp(t),n_E^\mp(t)\right)=\left(-\dot{f}(t)\mp\frac{f^2(t)n(t)}{2},\mp \dot{f}(t)-\frac{f^2(t)n(t)}{2}+n(t)\right).$$
\enD

We define the notion of a horocyclic parallel and examine its connection with horocyclic evolute.

\begin{definition}\label{def3.11} \rm
Let $(\bgamma,\bnu):I\rightarrow\Delta_1$ be a spacelike Legendre curve with curvature $(m,n)$. 
The \textit{horocyclic parallel} of $\bgamma$ (or, \textit{horocyclic parallel}  of $(\bgamma,\bnu)$) is   defined as
$$
P^\pm (\bgamma,\lambda_\pm)(t)=\bgamma(t)+\lambda _\pm(t)\bnu(t)+\frac{\lambda_\pm^2(t)}{2}(\bgamma(t)\pm\bmu(t)),
$$
where 
$ \lambda_\pm(t)  $ is a solution of the Bernoulli equation \begin{equation}\label{Bernoulli}
\frac{{\rm d}\lambda _\pm}{{\rm d}t}(t)=\mp \left(\frac{ \lambda _\pm^2(t) n(t)}{2}+\lambda _\pm(t)m(t)\right).
\end{equation}
\end{definition}

\begin{remark}\label{re3.12}\rm The solutions of   equation  \eqref{Bernoulli} are
\begin{equation}\label{lambda} 
 \lambda_\pm(t)=0 \text{ or } 	\lambda_\pm(t)=\frac{2e^{\int{\mp m(t) \mathrm{d}t}}}{\pm \int{n(t)e^{\int{\mp m(t) \mathrm{d}t}}}\mathrm{d}t+c}  \end{equation}for all $t\in I$, where $c$ is a constant. See also Remark \ref{Bernoulli1}.
\end{remark}

The horocyclic parallel of $\bgamma$ is a normal envelope of the one-parameter family of normal horocycles of $\bgamma$.

\begin{proposition} 
Let $(\bgamma,\bnu ):I\rightarrow\Delta_1$ be a spacelike Legendre curve  with curvature $(m,n)$ 
and $(\bx_N^\pm,{\bnu}_N^\pm):\R\times I\rightarrow \Delta_1$ is a one-parameter family of spacelike Legendre curves defined by Proposition \ref{p0}. 
Then the horocyclic parallel of $\bgamma$,	$$
P^\pm(\bgamma,\lambda_\pm)(t)=\bx_N^\pm\circ\be_N ^\pm[\pi/2] (t) =\bE_N ^\pm[\pi/2](t)
$$ is a $\pi/2$-enveloid  $($normal envelope$)$ of $\bx_N^\pm$, where  $\be_N^\pm[\pi/2]:I\rightarrow \R\times I $,
$ \be_N^\pm[\pi/2] (t)=(\lambda_\pm(t),t)$ is a pre-$\pi/2$-envelope of $\bx_N^\pm$ and $\lambda_\pm(t)$   is defined as equation \eqref{lambda}.
\end{proposition}
\demo Since $\lambda_\pm(t)$ is a solution of equation \eqref{Bernoulli}, the curvature of $(\bx_N^\pm,{\bnu}_N^\pm)$ satisfies 
$$\dot{\lambda}_\pm(t)m_N^\pm(\lambda_\pm(t),t)\sin({\pi}/{2})+M_N^\pm(\lambda_\pm(t),t)\sin({\pi}/{2})-L_N^\pm(\lambda_\pm(t),t)\cos({\pi}/{2})=0 $$ for all $t\in I$.
By Lemma \ref{enveloidtheorem}, we obtain $ \be_N^\pm[\pi/2] (t)=(\lambda_\pm(t),t)$ is a pre-$\pi/2$-envelope of $\bx_N^\pm$. Moreover, 
$$\begin{aligned}
P^\pm(\bgamma,\lambda_\pm)(t)&=\bgamma(t)+\lambda_\pm(t)\bnu(t)+\frac{\lambda_\pm^2(t)}{2}(\bgamma(t)\pm\bmu(t))
\\&= \bx_N^\pm (\lambda_\pm(t),t)=\bx_N^\pm\circ\be_N^\pm[\pi/2] (t) =\bE_N^\pm[\pi/2](t)
\end{aligned}$$ is a normal envelope of $\bx_N^\pm$.
\enD
\begin{proposition} \label{prop3.12}
If  $P^\pm(\bgamma,\lambda_\pm)$ is a  horocyclic parallel of $\bgamma$, then $(P^\pm(\bgamma,\lambda_\pm),\bnu_P^\pm):I\rightarrow\Delta_1$ is a spacelike Legendre curve with curvature
$\left(m_P^\pm(t),n_P^\pm(t)\right)=\left(\lambda _\pm(t)n(t)+m(t),n(t)\right),$
where $\bnu_P^\pm(t)=\lambda _\pm(t)\bgamma(t)+\bnu(t)\pm \lambda _\pm(t)\bmu(t).$ 
\end{proposition}
\demo  Take $\bnu_P^\pm(t)=\lambda _\pm(t)\bgamma(t)+\bnu(t)\pm \lambda _\pm(t)\bmu(t).$ 
By a direct calculation, we have $$ \dot{P}^\pm(\bgamma,{\lambda_\pm})(t)=(\lambda_\pm(t)n(t)+m(t)) \left(\mp( {\lambda _\pm^2(t)}/{2})\bgamma(t)\mp\lambda _\pm(t)\bnu(t)+\left(1- {\lambda _\pm^2(t)}/{2}\right)\bmu(t)\right).$$ It follows that   $ \langle P^\pm(\bgamma,\lambda_\pm)(t),\bnu_P^\pm(t)\rangle=\langle \dot{P}^\pm(\bgamma,\lambda_\pm)(t),\bnu_P^\pm(t)\rangle=0$  holds for all $t\in I$, which implies that $(P^\pm(\bgamma,\lambda_\pm),\bnu_P^\pm)$ is a spacelike Legendre curve.
Define $$\begin{aligned}
\bmu_P^\pm(t) =P^\pm(\bgamma,\lambda_\pm)(t) \wedge\bnu_P^\pm(t) 
 =\mp ({\lambda _\pm^2(t)}/{2})\bgamma(t)\mp\lambda _\pm(t)\bnu(t)+\left(1- {\lambda _\pm^2(t)}/{2}\right)\bmu(t).
\end{aligned}$$
Thus, $\{P^\pm(\bgamma,\lambda_\pm)(t),\bnu_P^\pm(t),\bmu_P^\pm(t)\}$ is a moving frame along $ P^\pm(\bgamma,\lambda_\pm)(t) $. Since $ \dot{P}^\pm(\bgamma,{\lambda_\pm})(t)=(\lambda_\pm(t)n(t)+m(t))\bmu_P^\pm(t)$  and $$ \dot{\bnu}_P^\pm(t)=n(t) \left(\mp( {\lambda _\pm^2(t)}/{2})\bgamma(t)\mp\lambda _\pm(t)\bnu(t)+\left(1- {\lambda _\pm^2(t)}/{2}\right)\bmu(t)\right)=n(t)\bmu_P^\pm(t),$$
the curvature of   $(P^\pm(\bgamma,\lambda_\pm),\bnu_P^\pm)$ is given by $\left(\lambda _\pm(t)n(t)+m(t),n(t)\right)$.
\enD
\begin{lemma}\label{lemma3.10}
Let $(\bgamma,\bnu ):I\rightarrow\Delta_1$ be a spacelike Legendre curve  with curvature $(m,n)$.	 Then the following assertions hold.  
\begin{enumerate}
\item[$(1)$] $m(t)+n(t)=0$ or  $m(t)-n(t)=0$ holds for all $t\in I$ if and only if there exists a lightlike constant vector $\ell$ such that $\bnu(t)\in \{\bx\in S_1^2\mid\langle\bx,\ell\rangle=1\}$ for all $t\in I$.
\item[$(2)$] Suppose that  the set of regular points of $\bgamma$ is   dense in $I$. If $\bgamma$ is a part of  horocycle, then  $m(t)+n(t)=0$ or  $m(t)-n(t)=0$ holds for all $t\in I$, and there exists a lightlike constant vector $  \ell$ such that $\bnu(t)\in \{\bx\in S_1^2\mid\langle\bx,\ell\rangle=1\}$  for all $t\in I$. 
\item[$(3)$] If either of the following two conditions is satisfied: {\rm(i)} $m(t)+n(t)=0$ or  $m(t)-n(t)=0$ holds for all $t\in I$. {\rm(ii)} There exists a lightlike constant vector $  \ell$ such that $\bnu(t)\in \{\bx\in S_1^2\mid\langle\bx,\ell\rangle=1\}$  for all $t\in I$, then $\bgamma$ is a   part of  horocycle.
\end{enumerate}
\end{lemma}
\demo {\rm (1)} 
If $m(t)+n(t)=0$ or  $m(t)-n(t)=0$ holds for all $t\in I$, then  $\dot{\bgamma}(t)=m(t)\bmu(t)=-n(t)\bmu(t)=-\dot{\bnu}(t)$ or $\dot{\bgamma}(t)=\dot{\bnu}(t)$ for all $t\in I$. 
Consequently, there exists a lightlike constant vector $\ell$ such that $\bgamma(t)=-\bnu(t)+\ell$ or $\bgamma(t)=\bnu(t)+\ell$ for all $t\in I$.  
It follows that $\left\langle \bnu(t),\ell \right\rangle=1$ or $\left\langle \bnu(t),-\ell \right\rangle=1$ for all $t\in I$.
\par
Conversely, suppose  there exists a lightlike constant vector $\ell$ such that $\left\langle \bnu(t),\ell \right\rangle=1$ for all $t\in I$.  
Differentiating this equality yields $\left\langle n(t)\bmu(t),\ell \right\rangle= n(t)\left\langle\bmu(t),\ell \right\rangle=0$. Since  the set $\mathrm{Reg}(\bnu)=\left\{t\in I\mid n(t)\neq0\right\}$  is  dense in $I$, we have  $\left\langle\bmu(t),\ell \right\rangle =0$ for all $t \in I$. Hence $\ell=\bgamma(t)+\bnu(t)$ or $\ell=-\bgamma(t)+\bnu(t)$ for all $t\in I$, so that $m(t)+n(t)=0$ or  $m(t)-n(t)=0$ holds for all $t\in I$.
\par
{\rm (2)} 
If $\bgamma$ is a part of  horocycle, there exists a   lightlike constant vector $\ell$ such that $\left\langle \bgamma(t),\ell \right\rangle=-1$ for all $t\in I$. Then $\left\langle m(t)\bmu(t),\ell \right\rangle= m(t)\left\langle\bmu(t),\ell \right\rangle=0$ for all $t\in I$. Since  the set of regular points of $\bgamma$ is  dense in $I$, we have  $\left\langle\bmu(t),\ell \right\rangle =0$ for all $t \in I$. 
Consequently, $\ell=\bgamma(t)+\bnu(t)$ or $\ell=\bgamma(t)-\bnu(t)$ for all $t\in I$,    which implies that  $m(t)+n(t)=0$ or $m(t)-n(t)=0$ for all $t\in I$. Moreover, $\ell$ satisifies  $\left\langle \bnu(t),\ell \right\rangle=1$ or $\left\langle \bnu(t),-\ell \right\rangle=1$ for all $t\in I$.
\par
{\rm (3)} If $m(t)+n(t)=0$ or  $m(t)-n(t)=0$ holds for all $t\in I$, the proof of (1) shows that
there exists a lightlike constant vector $\ell$ with  $\left\langle \bgamma(t),\ell \right\rangle=-1$  for all $t\in I$, which means that $\bgamma$ is a part of horocycle.
By (1), we see that conditions (i) and (ii) are equivalent. 
Therefore, if condition (ii) is satisfied, $\bgamma$ is also a part of horocycle.
\enD

\begin{proposition}\label{prop3.8}
Let $(\bgamma,\bnu ) :I\rightarrow\Delta_1$ be a spacelike Legendre curve   with curvature $(m,n)$.  $\left\{\bx_0,\bx_1,\bx_2\right\}$ is a pseudo-orthonormal basis of $\R_1^3$, where $\bx_0\in H^2$, $\bx_1,\bx_2\in S_1^2$ and $\bx_2=\bx_0\wedge\bx_1$.   Suppose that there exists a smooth function $f:I\rightarrow\R  $ such that  $m(t)+f(t)n(t)=0$ for all $t\in I$. Then the following assertions hold.
\begin{enumerate}
\item[$(1)$] ${Ev}^\pm(\bgamma) $ is a point if and only if $(\bgamma,\bnu)$ is congruent to a horocyclic parallel of $(\overline{\bgamma},\overline{\bnu} ):I\rightarrow\Delta_1$, where $\overline{\bgamma}(t)=\bx_0$, $\overline{\bnu}(t)=\cos\left(\int n(t) \mathrm{d}t\right)\bx_1+\sin\left(\int n(t) \mathrm{d}t\right)\bx_2$.

\item[$(2)$] Suppose that the set of regular points of ${Ev}^\pm(\bgamma)$ is  dense  in $I$.  If ${Ev}^\pm(\bgamma)$ is a part of  horocycle, then $(\bgamma,\bnu)$ is congruent to a horocyclic parallel of $(\widetilde{\bgamma},\widetilde{\bnu} ):I\rightarrow\Delta_1$, where \begin{equation}\label{eq4}\begin{aligned}
\widetilde{\bgamma}(t)&=\bx_0+\int n(t) \mathrm{d}t\;\bx_1+\frac{\left(\int n(t) \mathrm{d}t\right)^2}{2}\left(\bx_0+\bx_2\right),\\\widetilde{\bnu}(t)&=-\bx_2+\int n(t) \mathrm{d}t\;\bx_1+\frac{\left(\int n(t) \mathrm{d}t\right)^2}{2}\left(\bx_0+\bx_2\right).	\end{aligned}\end{equation} 

\item[$(3)$] Let $(\widetilde{\bgamma},\widetilde{\bnu}):I\rightarrow\Delta_1$ be given by equation \eqref{eq4}. If $(\bgamma,\bnu)$ is congruent to a horocyclic parallel of $(\widetilde{\bgamma},\widetilde{\bnu} )$,  then ${Ev}^\pm(\bgamma) $ is a part of  horocycle.
 \end{enumerate}
\end{proposition}
\demo {\rm(1)} 
Since $m_E^\pm(t)=-\dot{f}(t)\pm {f^2(t)n(t)}/{2}$,   ${Ev}^\pm(\bgamma) $ is a point if and only if  $ \dot{f}(t)=\pm {f^2(t)n(t)}/{2}$   for all $t\in I$. 
The solutions of the above equation are $f(t)=0$   or $f(t)={2}/{\left(-c\mp\int{n(t)}\mathrm{d}t\right)} $ for all $t\in I$, where $c $ is a constant. 
By $f(t)n(t)+m(t)=0$ for all $t\in I$ and the assumption that $\mathrm{Reg}(\bnu)=\left\{t\in I\mid n(t)\neq0\right\}$ is dense in $ I $,  we  obtain  the above result  is equivalent to $ m(t)=0 $ or $		m(t)=   {2n(t)}/{\left(\pm\int{n(t)}\mathrm{d}t+c\right)} $ for all $t\in I$. That is, the curvature of $(\bgamma,\bnu)$ is $$(0,n(t)) \text{ or }  \left(\frac{2n(t)}{\pm\int{n(t)}\mathrm{d}t+c},n(t)\right).$$ A direct calculation shows that the curvature of $(\overline{\bgamma},\overline{\bnu} ) $    is $(\overline{m}(t),\overline{n}(t))=(0,n(t))$, then  the solutions of equation $$\frac{{\rm d}\overline{\lambda} _\pm}{{\rm d}t}(t)=\mp \left(\frac{ \overline{\lambda} _\pm^2(t) \overline{n}(t)}{2}+\overline{\lambda} _\pm(t)\overline{m}(t)\right)$$ are $\overline{\lambda} _\pm(t)=0$   or $\overline{\lambda} _\pm(t)={2}/{\left(\pm\int{n(t)}\mathrm{d}t+c\right)} $ for all $t\in I$ by equation \eqref{lambda}. It follows by Proposition \ref{prop3.12} that the curvature   of the horocyclic parallel of $(\overline{\bgamma},\overline{\bnu} )$  satisfies $$ (\overline{m}^\pm_P(t),\overline{n}^\pm_P(t)) =(0,n(t))\text{ or }  (\overline{m}^\pm_P(t),\overline{n}^\pm_P(t)) =\left(\frac{2n(t)}{\pm\int{n(t)}\mathrm{d}t+c},n(t)\right).$$
Hence, there exists a $\overline{\lambda}_\pm(t)$ such that $(m(t),n(t))=(\overline{m}^\pm_P(t),\overline{n}^\pm_P(t))$ for all $t\in I$, which implies  $(\bgamma,\bnu)$ is congruent to a horocyclic parallel of $(\overline{\bgamma},\overline{\bnu} ) $.
\par
{\rm (2)}   If ${Ev}^\pm(\bgamma) $ is a  part of  horocycle, then by Lemma \ref{lemma3.10} and the assumption that the set of regular points of ${Ev}^\pm(\bgamma)$ is dense in $I$, we have   $m_E^\pm(t)\pm n_E^\pm(t)=\pm n(t)=0$ or $m_E^\pm(t)\mp n_E^\pm(t)=-2\dot{f}(t)\pm f^2(t)n(t)\mp n(t)=0$ for all $t\in I$.  Since $ \mathrm{Reg}(\bnu) $ is a dense subset of 
$ I $,  the first case cannot occur.  Thus, we obtain  $\dot{f} (t)=\pm{(f^2(t)-1)}n(t)/{2},$
the solutions of which are $f(t)=-1$ or $f(t)={\left(1+ce^{\int{\pm n(t)}\mathrm{d}t}\right)}/{\left(1-ce^{\int{\pm n(t)}\mathrm{d}t}\right)} $  for all $t\in I$. If follows that $m(t)-n(t)=0$ or  $m(t)= - {\left(1+ce^{\int{\pm n(t)}\mathrm{d}t}\right)n(t)}/{\left(1-ce^{\int{\pm n(t)}\mathrm{d}t}\right)}$ for all $t\in I$.   
That is, the curvature of $(\bgamma,\bnu)$ is $$ (n(t),n(t))\text{ or }  \left( -\frac{\left(1+ce^{\int{\pm n(t)}\mathrm{d}t}\right)n(t)}{1-ce^{\int{\pm n(t)}\mathrm{d}t}},n(t)\right).$$ 
By calculating, the curvature of $(\widetilde{\bgamma},\widetilde{\bnu} ) $    is $(\widetilde{m}(t),\widetilde{n}(t))=(n(t),n(t))$. Similar discussion yields that the curvature   of the horocyclic parallel of $\widetilde{\bgamma}$ is 
$$(\widetilde{m}^\pm_P(t),\widetilde{n}^\pm_P(t))= (n(t),n(t)) \text{ or }   (\widetilde{m}^\pm_P(t),\widetilde{n}^\pm_P(t))=\left(-\frac{\left(1+ce^{\int{\pm n(t)}\mathrm{d}t}\right)n(t)}{1-ce^{\int{\pm n(t)}\mathrm{d}t}}, n(t)\right).$$
Therefore, $(\bgamma,\bnu)$ is congruent to a horocyclic parallel of $(\widetilde{\bgamma},\widetilde{\bnu} ) $.
\par
{\rm (3)} The proof is the reverse process of that in {\rm (2)}.
\enD

\begin{proposition}\label{prop3.15} 
Suppose that there exists a smooth function $f:I\rightarrow\R  $ such that $m(t)+f(t)n(t)=0$ for all $t\in I$. Then the following assertions hold.   
\begin{enumerate}
\item[$(1)$] For any $t_0\in I$, there exists a horocyclic parallel $ P^\pm(\bgamma,\lambda_\pm) $ of $\bgamma$ such that $Ev^\pm(\bgamma)(t_0)= P^\pm(\bgamma,\lambda_\pm) (t_0)$ and $m_P^\pm(t_0)=0$. That is, the image of  $Ev^\pm(\bgamma)$   is   contained in the set of singular values of horocyclic parallels of $\bgamma$.
\par
\item[$(2)$] Let $ P^\pm(\bgamma,\lambda_\pm) $ be a horocyclic parallel satisfying $Ev^\pm(\bgamma)(t_0)= P^\pm(\bgamma,\lambda_\pm) (t_0)$, $m_P^\pm(t_0)=0$ as in {\rm(1)} and suppose that $n(t_0)\neq0$. Then $t_0$ is a regular point of $Ev^\pm(\bgamma)$   if and only if $t_0$ is a $(2,3)$-cusp of  $P^\pm(\bgamma,\lambda_\pm)$.
\par
\item[$(3)$]  Let $ P^\pm(\bgamma,\lambda_\pm) $ be a horocyclic parallel satisfying $Ev^\pm(\bgamma)(t_0)= P^\pm(\bgamma,\lambda_\pm) (t_0)$ and $m_P^\pm(t_0)=0$ as in {\rm(1)}. Then   $t_0$ is a $(2,3)$-cusp of $Ev^\pm(\bgamma)$  if and only if $t_0$ is a $(3,4)$-cusp of    $P^\pm(\bgamma,\lambda_\pm)$. 
\end{enumerate}
\end{proposition}
\demo {\rm(1)} For any $t_0\in I$, there exists a   solution  $ \lambda_\pm(t) $  of equation \eqref{Bernoulli} satisfying $ f(t_0)=\lambda _\pm(t_0)$. Thus, $Ev^\pm(\bgamma)(t_0)=P^\pm(\bgamma,\lambda_\pm)(t_0)$. Moreover,   $m_P^\pm(t_0)=\lambda _\pm(t_0) n(t_0)+m(t_0)=f(t_0) n(t_0)+m(t_0)=0$, so that $t_0$ is a singular point of  $ P^\pm(\bgamma,\lambda_\pm) $. 
\par
{\rm(2)}  Since $\dot{f}(t)n(t)+f(t)\dot{n}(t)+\dot{m}(t)=0$ holds for all $t\in I$ and $\lambda_\pm(t)$ is a solution of equation \eqref{Bernoulli},   we have
\begin{equation}\label{eq5}
\begin{aligned}
\dot{m}_P^\pm(t_0)&=\dot{\lambda}_\pm(t_0)n(t_0)+\lambda_\pm(t_0)\dot{n}(t_0)+\dot{m}(t_0)\\&=\mp\left(\frac{ \lambda _\pm^2(t_0) n(t_0)}{2}+\lambda _\pm(t_0)m(t_0)\right)n(t_0)+\lambda_\pm(t_0)\dot{n}(t_0)+\dot{m}(t_0)\\&=\pm\frac{f^2(t_0)n^2(t_0)}{2}-\dot{f}(t_0)n(t_0) =n(t_0)m_E^\pm(t_0).
\end{aligned}\end{equation}
Therefore, under the assumption $n(t_0)\neq0$, Corollary \ref{cor2.6} implies that $t_0$ is a regular point of $Ev^\pm(\bgamma)$ if and only if $t_0$ is a $(2,3)$-cusp of $P^\pm(\bgamma,\lambda_\pm)$.
\par
{\rm(3)} Direct calculating yields    \begin{equation}\label{eq6}\begin{aligned}
\ddot{m}_P^\pm(t_0)&=-\ddot{f}(t_0)n(t_0)\pm f(t_0)\dot{f}(t_0)n^2(t_0)\pm\frac{3f^2(t_0)\dot{n}(t_0)n(t_0)}{2}-2\dot{f}(t_0)\dot{n}(t_0)\\&=n(t_0) \dot{m}_E^\pm (t_0)+2\dot{n}(t_0)m_E^\pm(t_0).\end{aligned}\end{equation}
By Corollary \ref{cor2.6}, $t_0$ is a $(2,3)$-cusp of $Ev^\pm(\bgamma)$  if and only if  $m_E^\pm(t_0)=0$, $n_E^\pm(t_0)\neq0$ and $\dot{m}_E^\pm(t_0)\neq0$. According to equations \eqref{eq5}, \eqref{eq6} and Proposition \ref{prop3.1},  this is equivalent to   $n^\pm_P(t_0)=n(t_0)\neq0$, $\dot{m}_P^\pm(t_0)=0$  and $\ddot{m}_P^\pm(t_0)\neq0$. Since  $m_P^\pm(t_0)=0$, we obtain that   $t_0$ is a $(3,4)$-cusp of    $P^\pm(\bgamma,\lambda_\pm)$. 
\enD

\begin{proposition}\label{prop3.17}
Let $P^\pm(\bgamma,\lambda_\pm)$ be a horocyclic parallel of $\bgamma$.
Suppose that there exists a smooth function $f:I\rightarrow\R  $ such that $m(t)+f(t)n(t)=0$ for all $t\in I$. Then the horocyclic evolute of $P^\pm(\bgamma,\lambda_\pm) $ exists and
denote the horocyclic evolute of $(P^\pm(\bgamma,\lambda_\pm), \bnu_P^\pm)$ $($respectively, $(P^\pm(\bgamma,\lambda_\pm),-\bnu_P^\pm)$$)$ by  $Ev^\pm(P ^\pm(\bgamma,\lambda_\pm))$ $($respectively, $\overline{Ev}^\pm(P^\pm(\bgamma,\lambda_\pm))$$)$.     Furthermore, we have the following.
\begin{enumerate}
\item[$(1)$]  $Ev^+(P^+(\bgamma,\lambda_+))(t)=Ev^+(\bgamma)(t)$,    $Ev^-(P^-(\bgamma,\lambda_-))(t)= Ev^-(\bgamma)(t)$.
\item[$(2)$]  $\overline{Ev}^+(P^-(\bgamma,\lambda_-))(t)=Ev^-(\bgamma)(t) $,   $\overline{Ev}^-(P^+(\bgamma,\lambda_+))(t)=Ev^+(\bgamma)(t) $.
\end{enumerate} 
\end{proposition}
\demo {\rm(1)} Given that $f(t)n(t)+m(t)=0$ for all $t\in I$ and   $\mathrm{Reg}(\bnu)=\left\{t\in I\mid n(t)\neq0\right\}$ is   dense in $ I $,   Proposition \ref{prop3.12}  indicates the existence of a unique smooth function $f_P^\pm(t)=f(t)-\lambda _\pm(t)$ such that $f_P^\pm(t)n_P^\pm(t)+m_P^\pm(t)=0$ for all $t\in I$. Thus, the horocyclic evolute $Ev^+(P^+(\bgamma,\lambda_+))$ of $ (P^+(\bgamma,\lambda_+),\bnu_P^+) $ is $$\begin{aligned}
Ev^+(P^+(\bgamma,\lambda_+))(t)=&~P^+(\bgamma,\lambda_+)(t)+f_P^+(t)\bnu_P^+(t)+\frac{(f_P^+)^2(t)}{2}\left(P^+(\bgamma,\lambda_+)(t)+\bmu_P^+(t)\right)\\
=&~\bgamma(t)+\lambda _+(t)\bnu(t)+\frac{\lambda _+^2(t)}{2}(\bgamma(t)+\bmu(t))\\
&+(f(t)-\lambda _+(t))(\lambda _+(t)\bgamma(t)+\bnu(t)+\lambda _+(t)\bmu(t))\\&+\frac{(f(t)-\lambda _+(t))^2}{2}(\bgamma(t)+\bmu(t))\\
=&~\bgamma(t)+f(t)\bnu(t)+\frac{f^2(t)}{2}(\bgamma(t)+\bmu(t))=Ev^+(\bgamma)(t).
\end{aligned}$$
The equality $Ev^-(P^-(\bgamma,\lambda_-))(t)= Ev^-(\bgamma)(t)$   follows analogously.
\par
{\rm(2)} Note that $(P^\pm(\bgamma,\lambda_\pm),-\bnu_P^\pm)$ is a spacelike Legendre curve with curvature $(-m_P^\pm,n_P^\pm)$ and a moving frame $\{P^\pm(\bgamma,\lambda_\pm),-\bnu_P^\pm,-\bmu_P^\pm\}$. Applying Proposition \ref{prop3.7} and the results above, we obtain 
$
\overline{Ev}^-(P^+(\bgamma,\lambda_+))(t)= Ev^+(P^+(\bgamma,\lambda_+))(t)= Ev^+(\bgamma)(t) 
$
and
$ 
\overline{Ev}^+(P^-(\bgamma,\lambda_-))(t)= Ev^-(P^-(\bgamma,\lambda_-))(t)= Ev^-(\bgamma)(t). 
$
\enD
\section{Horocyclic involutes of spacelike frontals in $H^2$}\label{S4}

\begin{definition}\rm\label{def4.1}
Let $(\bgamma,\bnu ):I\rightarrow\Delta_1$ be a spacelike Legendre curve  with curvature $(m,n)$. The \textit{horocyclic involute} of $\bgamma$ (or, \textit{horocyclic involute}  of $(\bgamma,\bnu)$) is defined as 
$$
Inv^\pm(\bgamma,s_\pm)(t)=\bgamma(t)+s_\pm(t)\bmu(t)+\frac{s_\pm^2(t)}{2}(\bgamma(t)\pm\bnu(t)),
$$
where $s_\pm(t)$ is a solution of the Riccati equation
\begin{equation}\label{ODE}
\frac{{\rm d}s_\pm}{{\rm d}t}(t)= \frac{m(t)\pm n(t)}{2}s_\pm^2(t)-m(t).\end{equation}
\end{definition}

\begin{remark}{\rm
Once a particular solution of equation \eqref{ODE} is known, its general solution can be derived by the theory of Riccati equations.
}
\end{remark}

\begin{remark}\label{Bernoulli1}
\rm The equations \eqref{Bernoulli} and \eqref{ODE}  each admit  a unique local solution. That is, for any $t_0\in I$, there exists an interval $\overline{I} \subset I$ containing $t_0$
on which a solution exists and is unique. If a global solution does not exist, we examine the behavior of curves on such an  interval   $\overline{I}$.
\end{remark}

\begin{remark}{\rm
The \textit{horocyclic involute} of a regular curve $\bgamma:I\rightarrow H^2$ with geodesic curvature $\kappa_g$ is defined as 
$$
\mathcal{I}nv^\pm(\bgamma,s_\pm)(t)=\bgamma(t)+s_\pm(t)\bt(t)+\frac{s_\pm^2(t)}{2}(\bgamma(t)\mp\be(t)),
$$
where $s_\pm(t)$ is a solution of the Riccati equation
\begin{equation}\notag
\frac{{\rm d}s_\pm}{{\rm d}t}(t)= \frac{|\dot{\bgamma}(t)|(1\pm\kappa_g(t)) }{2}s_\pm^2(t)-|\dot{\bgamma}(t)|.
\end{equation}
 
}
\end{remark}

\begin{proposition}
Let $(\bgamma,\bnu ):I\rightarrow\Delta_1$ be a spacelike Legendre curve  with curvature $(m,n)$. Suppose that $t:\widetilde{I}\rightarrow I$ is a positive change of parameter, that is, $ t $ is surjective and has a positive derivative for all $u \in \widetilde{I}$. 
Then the \textit{horocyclic involute} of $\bgamma$ is  independent of the parameter change $t$. 
\end{proposition}
\demo
Let $(\widetilde{\bgamma}, \widetilde{\bnu}):\widetilde{I}\rightarrow\Delta_1$    be a spacelike Legendre curve defined by $(\widetilde{\bgamma}(u), \widetilde{\bnu}(u) )=(\bgamma(t(u)), \bnu(t(u)))$. According to Proposition 3.1 in \cite{CT2016},  its curvature is given by
$ 
(\widetilde{m}(u),\widetilde{n}(u))=(m(t(u))\dot{t}(u),n(t(u))\dot{t}(u)).
$ 
Multiplying  both sides of equation \eqref{ODE} by $\dot{t}(u)$ and setting $\widetilde{s}_\pm(u)=s_\pm(t(u))$  yields
\begin{equation}\notag
\begin{aligned}
\frac{{\rm d}s_\pm}{{\rm d}u}(t(u))= \frac{m(t(u))\dot{t}(u)\pm n(t(u))\dot{t}(u)}{2}s_\pm^2(t(u))-m(t(u))\dot{t}(u), 
\end{aligned}
\end{equation} 
that is,
$$
\frac{{\rm d}\widetilde{s}_\pm}{{\rm d}u}(u)=\frac{\widetilde{m }(u)\pm \widetilde{n}(u) }{2}\widetilde{s}^2_\pm(u)-\widetilde{m} (u).
$$
It follows that 
$$
\begin{aligned}
Inv^\pm(\bgamma,s_\pm)(t(u))&=\bgamma(t(u))+s_\pm(t(u))\bmu(t(u))+\frac{s_\pm^2(t(u))}{2}\left(\bgamma(t(u))\pm\bnu(t(u))\right)\\
&=\widetilde{\bgamma}(u)+\widetilde{s}_\pm(u)\widetilde{\bmu}(u)+\frac{\widetilde{s}^2_\pm(u)}{2}\left(\widetilde{\bgamma}(u)\pm\widetilde{\bnu}(u)\right)=Inv^\pm(\widetilde{\bgamma},\widetilde{s}_\pm)(u),
\end{aligned}
$$
where $\widetilde{\bmu}(u)=\bmu(t(u))$. This completes the proof.
\enD
Let $(\bgamma,\bnu):I\rightarrow \Delta_1$ be a spacelike Legendre curve with curvature $(m,n)$ and
take the frame
$ 
\left\{\bx_0(t),\bx_1(t), \bx_2(t)\right\}=\left\{\bgamma(t),\bmu (t),-\bnu(t) \right\}.
$ By Lemma \ref{lemma2.2}, we obtain
  the one-parameter family of tangent horocycles of $\bgamma$  
\begin{equation}\notag
\bx_T^\pm(s,t)=\bgamma(t)+s\bmu(t)+\frac{s^2}{2}(\bgamma(t)\pm\bnu(t)). 
\end{equation}

We can prove that  $\bx_T^\pm$ is a one-parameter family of spacelike frontals.
\begin{proposition}\label{p1}
Under the above notations, $(\bx_T^\pm,{\bnu}_T^\pm):\R\times I\rightarrow \Delta_1$ is a one-parameter family of spacelike Legendre curves with curvature 
$$
\begin{aligned}
&\left({m}_T^\pm(s,t),{n}_T^\pm(s,t),{L}_T^\pm(s,t),{M}_T^\pm(s,t),{N}_T^\pm(s,t)\right)\\=&\left(\mp1,\mp1,s(m(t)\pm n(t)),\mp m(t)\pm s^2(m(t)\pm n(t))/2,n(t)\pm s^2(m(t)\pm n(t))/2\right),	
\end{aligned}
$$
where 
${\bnu}_T^\pm(s,t)=\mp\bnu(t)+s\bmu(t)+({s^2}/{2})(\bgamma(t)\pm\bnu(t)).$
\end{proposition}
\demo Take ${\bnu}_T^\pm(s,t)=\mp\bnu(t)+s\bmu(t)+({s^2}/{2})(\bgamma(t)\pm\bnu(t)).$ A direct calculation yields $\bx_{Ts}^\pm(s,t)= \bmu(t)+s(\bgamma(t)\pm\bnu(t)). $
It follows that 
$ 
\langle  \bx_T^\pm(s,t),\bnu_T^\pm(s,t)\rangle=\langle   \bx_{Ts}^\pm(s,t),\bnu_T^\pm(s,t)\rangle=0
$   
holds for all $(s,t)\in \R\times I$, which implies that $( \bx_T^\pm,\bnu_T^\pm)$ is a    one-parameter family of spacelike Legendre curves. 
Define 
$\bmu_T^\pm(s,t)=\bx_T^\pm(s,t)\wedge\bnu_T^\pm(s,t)=-s\bnu(t)\mp\bmu(t)\mp s\bgamma(t).$ Then $\{\bx_T^\pm(s,t),\bnu_T^\pm(s,t),\bmu_T^\pm(s,t)\}$ is a  moving frame along $\bx_T^\pm(s,t)$.  By further calculations, we obtain 
$$\begin{aligned}\bx_{Ts}^\pm(s,t)&=\bnu_{Ts}^\pm(s,t)= \bmu(t)+s(\bgamma(t)\pm\bnu(t))=\mp\bmu_T^\pm(s,t),\\
\bx_{Tt}^\pm(s,t)&=sm(t)\bgamma(t)-sn(t)\bnu(t)+((1+s^2/2)m(t)\pm s^2n(t))\bmu(t) \\
&=s(m(t)\pm n(t))\bnu_T^\pm(s,t)+(\mp m(t)\pm s^2(m(t)\pm n(t))/2) \bmu_T^\pm(s,t),\\
\bnu_{Tt}^\pm(s,t)&=sm(t)\bgamma(t)-sn(t)\bnu(t)+(s^2m(t)/2\mp(1- s^2)n(t))\bmu(t)\\
&=s(m(t)\pm n(t))\bx_T^\pm(s,t)+(n(t)\pm s^2(m(t)\pm n(t))/2) \bmu_T^\pm(s,t).\end{aligned}$$
Therefore,  the curvature   of one-parameter family of spacelike Legendre curves $( \bx_T^\pm,\bnu_T^\pm)$ is given by $\left(\mp1,\mp1,s(m(t)\pm n(t)),\mp m(t)\pm s^2(m(t)\pm n(t))/2,n(t)\pm s^2(m(t)\pm n(t))/2\right)$.
\enD

\begin{proposition}\label{prop4.6}
Let $(\bgamma,\bnu ):I\rightarrow\Delta_1$ be a spacelike Legendre curve  with curvature $(m,n)$ 
and $(\bx_T^\pm,{\bnu}_T^\pm):\R\times I\rightarrow \Delta_1$ is a one-parameter family of spacelike Legendre curves defined by Proposition \ref{p1}. 
Then the horocyclic involute of $\bgamma$, $$
Inv^\pm(\bgamma,s_\pm)(t) =\bx_T^\pm\circ\be_T^\pm[\pi/2](t)=\bE_T^\pm[\pi/2](t)
$$ is a $\pi/2$-enveloid  $($normal envelope$)$ of $\bx_T^\pm$, where $\be_T^\pm[\pi/2]:I\rightarrow \R\times I $,
$ \be_T^\pm [\pi/2](t)=(s_\pm(t),t)$ is a pre-$\pi/2$-enveloid of $\bx_T^\pm$ and $s_\pm(t)$ is a solution of equation \eqref{ODE}.

\end{proposition}
\demo
Since $s_\pm(t)$   is a solution of equation \eqref{ODE},
the curvature of $(\bx_T^\pm,{\bnu}_T^\pm)$ satisfies
$$
\dot{s}_\pm(t){m}_T^\pm(s_\pm(t),t)\sin({\pi}/{2})+ {M}_T^\pm(s_\pm(t),t)\sin({\pi}/{2})-{L}_T^\pm(s_\pm(t),t)\cos({\pi}/{2})=0 
$$ for all $t\in I$.
Hence, by Lemma \ref{enveloidtheorem},   $\be_T^\pm[\pi/2](t)=(s_\pm(t),t)$ is a pre-$\pi/2$-enveloid of $ \bx_T^\pm $. 
Moreover,
\begin{equation}\notag
\begin{aligned}
Inv^\pm(\bgamma,s_\pm)(t)
&=\bgamma(t)+s_\pm(t)\bmu(t)+\frac{s_\pm^2(t)}{2}(\bgamma(t)\pm\bnu(t))\\ &={\bx}_T^\pm(s_\pm(t),t)  = \bx_T^\pm\circ\be_T^\pm[\pi/2](t)=  {\bE}_T^\pm[\pi/2](t)
\end{aligned}
\end{equation}
is a  normal envelope of $\bx_T^\pm$.
\enD
Given a solution $s_\pm(t)$ of equation \eqref{ODE}, consider the function
  $H_I^\pm:I\times H^2\rightarrow\R$,  $$H_I^\pm(t,\bx)=\left\langle \bgamma(t)\mp\bmu(t)+s_\pm(t)\left(\bmu(t)\mp\bgamma(t)-\bnu(t)\right)+\frac{s^2_\pm(t)}{2}(\bgamma(t)\pm\bnu(t)),\bx\right\rangle+1.$$   The following conclusion can be drawn.

\begin{proposition}\label{prop4.7}
Under the above notations, suppose that $( s_\pm(t)\mp1)(n(t)\mp m(t))\neq0$ for all $t\in I$, then   the image of horocyclic involute  ${Inv}^\pm(\bgamma,s_\pm) $  of $\bgamma$ is the discriminant set of $ H_I^\pm$.
\end{proposition}
\demo For any $(t, \bx)\in I\times H^2$, there exist $a,b,c\in\R$ such that $\bx=a\bgamma(t)+b
\bmu(t)+c\bnu(t)$.
Under the condition $H_I^\pm(t,\bx) =0 $
and  $\langle\bx,\bx\rangle=-a^2+b^2+c^2 =-1$, we  find that there exists $\lambda\in\R$ such that 
$$\begin{aligned}
a&=\left(1+\frac{\lambda^2}{2}\right)\left(1+\frac{s_\pm^2(t)}{2}\right)\mp\frac{\lambda s_\pm^2(t)}{2}\mp\frac{\lambda^2 s_\pm(t)}{2},\\
b&=\left(1+\frac{\lambda^2}{2}\right)s_\pm(t)\mp \lambda s_\pm(t)\mp\frac{\lambda^2}{2},\\
c&=\pm\left(1+\frac{\lambda^2}{2}\right)\frac{s_\pm^2(t)}{2}+\lambda\left(1-\frac{s_\pm^2(t)}{2}\right)-\frac{\lambda^2 s_\pm(t)}{2}.
\end{aligned}$$
Furthermore, differentiating $H_I^\pm(t,\bx)$ with respect to $t$ and substituting $\bx$ yields
$$\begin{aligned} \frac{\partial H_I^\pm}{\partial t}(t,\bx)&=\left\langle  -( s_\pm(t)\mp1)(n(t)\mp m(t))\left(\bnu(t)\mp s_\pm(t)\bmu(t)\mp\frac{s^2_\pm(t)}{2}(\bgamma(t)\pm\bnu(t))\right),\bx\right\rangle\\
&=-( s_\pm(t)\mp1)(n(t)\mp m(t))\lambda.\end{aligned}$$
By $\left({\partial H_I^\pm}/{\partial t}\right)(t,\bx)=0$, we obtain $\lambda=0$ under the condition $( s_\pm(t)\mp1)(n(t)\mp m(t))\neq0$ for all $t\in I$.  Hence, $$\bx=\left(1+\frac{s_\pm^2(t)}{2}\right)\bgamma(t)+s_\pm (t)\bmu(t)\pm\frac{s_\pm^2(t)}{2}\bnu(t)=Inv^\pm(\bgamma,s_\pm)(t).$$ This completes the proof of the proposition.
\enD

\begin{proposition}\label{prop4.8}
If ${Inv}^\pm(\bgamma,s_\pm)$ is a horocyclic involute of $\bgamma$, then  
$\left({Inv}^\pm(\bgamma,s_\pm), \bnu_I^\pm \right):I\rightarrow \Delta_1$ is a spacelike Legendre curve with curvature  
$$
\left(m^\pm_I(t), n^\pm_I(t)\right)= \left(- s_\pm(t)(n(t)\pm m(t)),n(t)\pm m(t)\right),
$$ 
where 
$ 
\bnu^\pm_I(t)=- s_\pm(t)\bgamma(t)-\bmu(t)\mp s_\pm(t)\bnu(t).
$ 
\end{proposition}
\demo Take $ 
\bnu^\pm_I(t)=- s_\pm(t)\bgamma(t)-\bmu(t)\mp s_\pm(t)\bnu(t).
$ Since $${\dot{I}nv}^\pm(\bgamma,s_\pm)(t)=-s_\pm(t)(n(t)\pm m(t))\left(\bnu(t)\mp s_\pm(t)\bmu(t)\mp ({s^2_\pm(t)}/{2})(\bgamma(t)\pm\bnu(t))\right),$$
 we have
$
\langle {Inv}^\pm(\bgamma,s_\pm)(t),\bnu^\pm_I(t)\rangle=\langle {\dot{I}nv}^\pm(\bgamma,s_\pm)(t),\bnu^\pm_I(t)\rangle=0
$
for all $t\in I$, which implies that  $\left({Inv}^\pm(\bgamma,s_\pm), \bnu_I^\pm \right) $ is a spacelike Legendre curve. 
Define
$$
\begin{aligned}
\bmu_I^\pm(t) ={Inv}^\pm(\bgamma,s_\pm)(t)\wedge \bnu_I^\pm(t) 
 =\bnu(t)\mp s_\pm(t)\bmu(t)\mp({s^2_\pm(t)}/{2})(\bgamma(t)\pm\bnu(t)).
\end{aligned}
$$ 
Then $\left\{{Inv}^\pm(\bgamma,s_\pm)(t), \bnu_I^\pm(t), \bmu_I^\pm(t)\right\}$ is a moving frame along ${Inv}^\pm(\bgamma,s_\pm)(t)$. A direct calculation   yields  $ {\dot{I}nv}^\pm(\bgamma,s_\pm)(t)=-s_\pm(t)(n(t)\pm m(t)) \bmu_I^\pm(t) $  and $$\dot{\bnu}_I^\pm(t)=(n(t)\pm m(t))\left(\bnu(t)\mp s_\pm(t)\bmu(t)\mp ({s^2_\pm(t)}/{2})(\bgamma(t)\pm\bnu(t))\right)=(n(t)\pm m(t)) \bmu_I^\pm(t).$$ Then the curvature of  $\left({Inv}^\pm(\bgamma,s_\pm), \bnu_I^\pm \right) $ is given by $\left(- s_\pm(t)(n(t)\pm m(t)),n(t)\pm m(t)\right)$.
\enD

Combining Corollary \ref{cor2.6} and Proposition \ref{prop4.8}, we have the following result.
\begin{corollary}\label{cor4.8}
Let $(\bgamma,\bnu ):I\rightarrow\Delta_1$ be a spacelike Legendre curve  with curvature $(m,n)$.	   ${Inv}^\pm(\bgamma,s_\pm) $ is a horocyclic involute of $ \bgamma$ and $t_0\in I$.
\begin{enumerate}
\item[$(1)$] ${Inv}^\pm(\bgamma,s_\pm)$ is singular at $t_0$ if and only if $ s_\pm(t_0)(n(t_0)\pm m(t_0))=0$.

\item[$(2)$] ${Inv}^\pm(\bgamma,s_\pm)$ has a $(2,3)$-cusp at $t_0$ if and only if $ s_\pm(t_0)=0$, $m(t_0)\neq0$, $  n(t_0)\pm m(t_0) \neq0$.

\item[$(3)$] ${Inv}^\pm(\bgamma,s_\pm)$ has a $(3,4)$-cusp at $t_0$ if and only if  $ s_\pm(t_0)=0$, $m(t_0)=0$, $n(t_0)\neq0$ and $\dot{m}(t_0)\neq0$.

\item[$(4)$] ${Inv}^\pm(\bgamma,s_\pm)$ has a $(2,5)$-cusp at $t_0$ if and only if $n(t_0)\pm m(t_0) =0$, $\dot{n}(t_0)\pm \dot{m}(t_0) \neq0$, $ s_\pm(t_0)\neq0$ and $m(t_0)\neq0$.

\item[$(5)$] ${Inv}^\pm(\bgamma,s_\pm)$ has a $(3,5)$-cusp at $t_0$ if and only if $ s_\pm(t_0)=n(t_0)\pm m(t_0) =0$, $ \dot{n}(t_0)\pm \dot{m}(t_0) \neq0$ and $m(t_0)\neq0$.
\end{enumerate}
\end{corollary}

\begin{corollary}\label{cor4.9}
Let $(\bgamma,\bnu ):I\rightarrow\Delta_1$ be a spacelike Legendre curve  with curvature $(m,n)$.	   ${Inv}^\pm(\bgamma,s_\pm) $ is a horocyclic involute of $ \bgamma$ and $t_0\in I$.
\begin{enumerate}
\item[$(1)$] If ${Inv}^\pm(\bgamma,s_\pm)$ intersects $\bgamma$ at $t_0$, then $t_0$ is a singular point of ${Inv}^\pm(\bgamma,s_\pm)$. 

\item[$(2)$] If $n(t)\pm m(t)\neq0$ for all $t\in I$, then ${Inv}^\pm(\bgamma,s_\pm) $ is a spacelike front and has no inflection point.

\item[$(3)$] If $t_0$ is a $(2,3)$-cusp, $(2,5)$-cusp or $(3,5)$-cusp of $Inv^\pm(\bgamma,s_\pm)$, then $t_0$ is a regular point of $\bgamma$.

\item[$(4)$] If  $t_0$ is a $(3,4)$-cusp of $Inv^\pm(\bgamma,s_\pm)$, then  $t_0$ is a $(2,3)$-cusp of $\bgamma$.

\item[$(5)$]  If $t_0$ is a   $(2,5)$-cusp or $(3,5)$-cusp of $ \bgamma $, then $t_0$ is a singular point of $Inv^\pm(\bgamma,s_\pm)$.
\end{enumerate}
\end{corollary}

Since $s_\pm(t)$ is a solution of equation \eqref{ODE}, then
$\overline{s}_\pm(t)=-s_\mp(t)$ is a solution of the Riccati equation 
\begin{equation}\label{eq9}
\frac{{\rm d}\overline{s}_\pm}{{\rm d}t}(t) = \frac{ - m(t)\pm n(t)}{2}\overline{s}_\pm^2(t)+(-m(t)).  \end{equation}
Denoting the horocyclic involute of $(\bgamma,-\bnu)$ by $\overline{Inv}^\pm(\bgamma,\overline{s}_\pm)$ in this paper, we obtain   the following proposition.
\begin{proposition}\label{prop4.10}
$(\overline{Inv}^\pm(\bgamma,\overline{s}_\pm),\overline{\bnu}_I^\pm):I\rightarrow\Delta_1$ is a spacelike Legendre curve with curvature $\left(- s_\mp(t)(n(t)\mp m(t)),n(t)\mp m(t)\right),$
where $\overline{Inv}^\pm(\bgamma,\overline{s}_\pm)=Inv^\mp(\bgamma,s_\mp)$
and $\overline{\bnu}_I^\pm=\bnu_I^\mp$.
\end{proposition}
\demo Note that $(\bgamma,-\bnu) $ is a spacelike Legendre curve with curvature $(-m,n)$ and a moving frame $\{\bgamma,-\bnu,-\bmu\}$.
Since   $\overline{s}_\pm(t)$ is a solution of equation \eqref{eq9}, we can get
 the horocyclic involute of $(\bgamma,-\bnu)$  is
$$
\begin{aligned}
\overline{Inv}^\pm(\bgamma,\overline{s}_\pm)(t)&=\bgamma(t)+\overline{s}_\pm(t)(-\bmu(t))+\frac{\overline{s}^2_\pm(t)}{2}(\bgamma(t)\mp\bnu(t))\\
&=\bgamma(t)+s_\mp(t)\bmu(t)+\frac{s_\mp^2(t)}{2}(\bgamma(t)\mp\bnu(t)) =Inv^\mp(\bgamma,s_\mp)(t).
\end{aligned}$$
Take $\overline{\bnu}_I^\pm=\bnu_I^\mp$. By Proposition \ref{prop4.8}, $(\overline{Inv}^\pm(\bgamma,\overline{s}_\pm),\overline{\bnu}_I^\pm)=(Inv^\mp(\bgamma,s_\mp),\bnu_I^\mp)$ is a spacelike Legendre curve with curvature  $\left(\overline{m}_I^\pm(t),\overline{n}_I^\pm(t)\right)=\left(m_I^\mp(t),n_I^\mp(t)\right)= \left(- s_\mp(t)(n(t)\mp m(t)),n(t)\mp m(t)\right).$
\enD

\begin{proposition}\label{prop4.13}
For any solutions  $s_\pm(t)$ and $\widehat{s}_\pm(t) $ of equation \eqref{ODE}, the horocyclic involutes \\ $Inv^\pm(\bgamma,s_\pm)$ and $Inv^\pm(\bgamma,\widehat{s}_\pm)$ of $\bgamma$ are
horocyclic parallels.
\end{proposition}
\demo
By Definition \ref{def4.1} and Proposition \ref{prop4.8}, we have
\begin{equation}\notag
\begin{aligned}
&Inv^\pm(\bgamma,\widehat{s}_\pm)(t)=~\bgamma(t)+\widehat{s}_\pm(t)\bmu(t)+\frac{\widehat{s}_\pm^2(t)}{2}(\bgamma(t)\pm\bnu(t))\\
=&~\bgamma(t)+s_\pm(t)\bmu(t)+\frac{s_\pm^2(t)}{2}(\bgamma(t)\pm\bnu(t))+(s_\pm(t)- \widehat{s}_\pm(t))\left(-s_\pm(t)\bgamma(t)-\bmu(t)\mp s_\pm(t)\bnu(t)\right)\\&~+\frac{(s_\pm(t) -\widehat{s}_\pm(t))^2}{2}\left(\bgamma(t)\pm
\bnu(t)\right)  \\
=&~Inv^\pm(\bgamma,s_\pm)(t)+(s_\pm(t)- \widehat{s}_\pm(t))\bnu^\pm_I(t)+\frac{(s_\pm(t) -\widehat{s}_\pm(t))^2}{2}\left(Inv^\pm(\bgamma,s_\pm)(t)\pm
\bmu_I^\pm(t)\right).
\end{aligned}
\end{equation}
Furthermore,  
\begin{equation}\notag
\begin{aligned}
\frac{{\rm d}(s_\pm(t)- \widehat{s}_\pm(t))}{{\rm d}t}&= \frac{m(t)\pm n(t)}{2}\left(s_\pm^2(t)-\widehat{s}^2_\pm(t)\right)\\
&=\mp\left(\frac{n(t)\pm m(t)}{2}(s_\pm(t)-\widehat{s}_\pm(t))^2-s_\pm(t)(n(t)\pm m(t))(s_\pm(t)-\widehat{s}_\pm(t))\right)\\
&=\mp\left(\frac{n^\pm_I(t)}{2}(s_\pm(t)-\widehat{s}_\pm(t))^2+m_I^\pm(t)(s_\pm(t)-\widehat{s}_\pm(t))\right).
\end{aligned}
\end{equation}
The solutions of the  above equation are $$s_\pm(t)- \widehat{s}_\pm(t)=0 
\text{ or } s_\pm(t)- \widehat{s}_\pm(t)= \frac{2e^{\int{\mp m^\pm_I(t) \mathrm{d}t}}}{\pm\int{n^\pm_I(t)e^{\int{\mp m^\pm_I(t) \mathrm{d}t}}}\mathrm{d}t+c}  $$   for all $t\in I$, where $c $ is a constant. 
By the definition of horocyclic parallels, $ Inv^\pm(\bgamma,\widehat{s}_\pm) $ is a horocyclic  parallel  of $ Inv^\pm(\bgamma,s_\pm) $. Similarly, $Inv^\pm(\bgamma,s_\pm)$ is also a horocyclic parallel of $Inv^\pm(\bgamma,\widehat{s}_\pm)$.
\enD

\begin{proposition}
Let $(\bgamma,\bnu ):I\rightarrow\Delta_1$ be a spacelike Legendre curve  with curvature $(m,n)$.	   ${Inv}^\pm(\bgamma,s_\pm) $ is a horocyclic involute of $ \bgamma$. 
\begin{enumerate}
\item[$(1)$] If $({Inv}^\pm(\bgamma,s_\pm),\bnu_I^\pm) $ is  real analytic and ${Inv}^\pm(\bgamma,s_\pm) $ is a point,    then   $\bgamma$ is a point or  a part of horocycle.   

\item[$(2)$] If $ \bgamma $ is a point, then   ${Inv}^\pm(\bgamma,s_\pm) $ is a  horocyclic parallel  of a point.

\item[$(3)$] Suppose that $(Inv^\pm(\bgamma,s_\pm),\bnu_I^\pm )$ is real analytic  and the set of regular points of $ Inv^\pm(\bgamma,s_\pm)$ is  dense in   $I$. If $ Inv^\pm(\bgamma,s_\pm)$ is a part of horocycle, then $ \bgamma $ is a part of horocycle.

\item[$(4)$] Suppose that the set of regular points of $  \bgamma  $ is  dense in $I$. If $ \bgamma $ is a  part of  horocycle, then  ${Inv}^\pm(\bgamma,s_\pm) $  is a point on a horocycle or a horocyclic parallel  of  a part of horocycle.
\end{enumerate}
\end{proposition}
\demo{\rm(1)} Since $m_I^\pm(t)=-s_\pm(t)(n(t)\pm m(t))$ is  real analytic and   ${Inv}^\pm(\bgamma,s_\pm) $ is a point,  we have $s_\pm(t)=0$ or $n(t)\pm m(t)=0$ for all $t\in I$.   If $s_\pm(t)=0$ for all $t\in I$, equation \eqref{ODE} gives $m(t)=0$ for all $t\in I$. Then  $\bgamma$ is a point. 
In the case when  $n(t)\pm m(t)=0$ for all $t\in I$,   
 $\bgamma$ is a part of horocycle   by Lemma \ref{lemma3.10}. 

{\rm(2)} If  $ \bgamma $ is a point,  then $m(t)=0$ for all $t\in I$.  Thus, the solutions of equation \eqref{ODE} are $$s_\pm(t)=0 \text{ or }  s_\pm(t)=\frac{2}{\mp\int{n(t)\mathrm{d}t+c}}$$  for all $t\in I$, where $c $ is a constant. It follows that  the horocyclic involute satisfying $s_\pm(t)=0$ for all $t\in I$ is a point and hence ${Inv}^\pm(\bgamma,s_\pm) $ is a horocyclic parallel  of a point by Proposition \ref{prop4.13}.

{\rm (3)}  Taking $Inv^+(\bgamma,s_+)  $  
as an example, the proof for $Inv^-(\bgamma,s_-)  $ is analogous. By Lemma \ref{lemma3.10}, the curvature of $(Inv^+(\bgamma,s_+),\bnu_I^+ )$ satisfies $m_I^+(t)+n_I^+(t)=(1-s_+(t))(n(t)+m(t))=0$ or  $m_I^+(t)-n_I^+(t)=-(1+s_+(t))(n(t)+m(t))=0$ for all $t\in I$. Since $(Inv^+(\bgamma,s_\pm),\bnu_I^+ )$ is real analytic, we obtain the following three cases: $1-s_+(t)=0$,  $1+s_+(t)=0$  or $n(t)+m(t)=0$ for all $t\in I$.  Substituting these into equation \eqref{ODE} yields $n(t)-m(t)=0$ or $n(t)+m(t)=0$  for all $t\in I$. Thus, we conclude that $ \bgamma  $   is a  part of  horocycle by Lemma  \ref{lemma3.10}. 

{\rm (4)}   If $ \bgamma  $ is a  part of  horocycle, then by Lemma  \ref{lemma3.10} and the assumption that the set of regular points of $\bgamma$ is  dense  in $I$, we have either  $m (t)+n (t)= 0$ or  $m (t)-n (t) =0$  for all $t\in I$.  
If $m (t)+n (t)= 0$ for all $t\in I$, we obtain $m_I^+(t)=n_I^+(t)=0$ for all $t\in I$. It follows that $Inv^+(\bgamma,s_+)$ is a point on a horocycle. 
If $m (t)- n (t) =0$ for all $t\in I$, the solutions of equation \eqref{ODE} are $$s_\pm(t)=-1   \text{ or  } s_\pm(t)=\frac{1+ce^{2\int{\pm n(t)}\mathrm{d}t}}{1-ce^{2\int{\pm n(t)}\mathrm{d}t}} $$ for all $t\in I$. 
For $c=0$, we have $s_+(t)=1$ and $m_I^+(t)+n_I^+(t) =0$ for all $t\in I$. Consequently, the corresponding horocyclic involute ${Inv}^+(\bgamma,1) $  is a part of  horocycle and thus ${Inv}^+(\bgamma,s_+) $ is a horocyclic parallel  of a part of horocycle by Proposition \ref{prop4.13}. The proof for $Inv^-(\bgamma,s_-)  $  follows similarly.
\enD

We give  relations between horocyclic evolutes and   involutes of $\bgamma$.
\begin{theorem}\label{th4.14}
$(1)$ We denote $\mathrm{Reg}(\bnu^\pm_I)=\left\{t\in I\mid n(t)\pm m(t)\neq0\right\}$ as the set of regular points of $\bnu_I^\pm$. If  $\mathrm{Reg}(\bnu^\pm_I)$ is  dense in $I$, then the horocyclic evolute of the horocyclic involute of  $\bgamma$ always exists and is unique.   Moreover,  $$
\begin{aligned}
Ev^+\left(Inv^+(\bgamma,s_+)\right)(t)&=Ev^-\left(Inv^-(\bgamma,s_-)\right)(t)=\bgamma(t),\\
\overline{Ev}^+\left(Inv^-(\bgamma,s_-)\right)(t)&=\overline{Ev}^-\left(Inv^+(\bgamma,s_+)\right)(t)=\bgamma(t).\\
\end{aligned}
$$
$(2)$ Suppose that there exists a smooth function $f:I\rightarrow\R  $ such that $m(t)+f(t)n(t)=0$ for all $t\in I$. Then
$$
\begin{aligned}
Inv^+\left(Ev^+(\bgamma),S^+_E\right)(t)&=\overline{Inv}^-\left(Ev^+(\bgamma),\overline{S}^-_E\right)(t)\\&=  \bgamma(t)+\left(f(t)-S^+_E(t)\right)\bnu(t)+\frac{\left(f(t)-S^+_E(t)\right)^2}{2}\left(\bgamma(t)+\bmu(t)\right),\\
Inv^-\left(Ev^-(\bgamma),S^-_E\right)(t)&= \overline{Inv}^+\left(Ev^-(\bgamma),\overline{S}^+_E\right)(t)\\&=  \bgamma(t)+\left(f(t)-S^-_E(t)\right)\bnu(t)+\frac{\left(f(t)-S^-_E(t)\right)^2}{2}\left(\bgamma(t)-\bmu(t)\right),
\end{aligned}
$$
where 
$$
S_E^\pm(t)=f(t) \text{ or }
S_E^\pm(t)=f(t)-\frac{2e^{\mp\int m(t)\mathrm{d}t}}{\pm\int n(t)e^{\mp\int m(t)\mathrm{d}t}\mathrm{d}t+c},\; \overline{S}_E^\pm(t)=-S_E^\mp(t) 
$$
and $c$ is a constant. 
Moreover, $Inv^+ (Ev^+(\bgamma),S^+_E)$, $Inv^- (Ev^-(\bgamma),S^-_E)$, $\overline{Inv}^+ (Ev^-(\bgamma),\overline{S}^+_E)$ and 
$\overline{Inv}^- (Ev^+(\bgamma),\overline{S}^-_E)$ are horocyclic parallels of $\bgamma$.	
\end{theorem}
\demo {\rm(1)} Since the set $\mathrm{Reg}(\bnu^\pm_I)$ is  dense in $I$,  Proposition \ref{prop4.8} implies the existence of a unique smooth function $f_I^\pm(t)=s_\pm(t)$ satisfying $f_I^\pm(t)n_I^\pm(t)+m_I^\pm(t)=0$ for all $t\in I$. 
Then
$$
\begin{aligned}
Ev^+\left(Inv^+(\bgamma,s_+)\right)(t)=&~Inv^+(\bgamma,s_+)(t)+f_I^+(t)\bnu_I^+(t)+\frac{(f_I^+)^2(t)}{2}\left(Inv^+(\bgamma,s_+)(t)+\bmu_I^+(t)\right)\\
=&~\bgamma(t)+s_+(t)\bmu(t)+\frac{s_+^2(t)}{2}(\bgamma(t)+\bnu(t))\\
&+s_+(t)\left(- s_+(t)\bgamma(t)-\bmu(t)- s_+(t)\bnu(t)\right) +\frac{s_+^2(t)}{2}(\bgamma(t)+ \bnu(t))\\
=&~\bgamma(t).
\end{aligned}
$$
Similarly, one shows that $Ev^-\left(Inv^-(\bgamma,s_-)\right)(t)=\bgamma(t)$.

Since  $\left({Inv}^\pm(\bgamma,s_\pm), -\bnu_I^\pm \right) $ is a spacelike Legendre curve with curvature  
$\left(-m^\pm_I, n^\pm_I \right)$ and a moving frame $\{{Inv}^\pm(\bgamma,s_\pm), -\bnu_I^\pm,-\bmu_I^\pm \}$. Applying Proposition \ref{prop3.7} together with the results above, we obtain 
\begin{align*}
\overline{Ev}^+\left(Inv^-(\bgamma,s_-)\right)(t) &= {Ev}^-\left(Inv^-(\bgamma,s_-)\right)(t)=\bgamma(t),\\
\overline{Ev}^-\left(Inv^+(\bgamma,s_+)\right)(t) &= {Ev}^+\left(Inv^+(\bgamma,s_+)\right)(t)=\bgamma(t).
\end{align*}

{\rm(2)} By Proposition \ref{prop3.1},
$\left({Ev}^\pm(\bgamma ), \bnu^\pm_E \right) $ is a spacelike Legendre curve with curvature 
$$
\left(m^\pm_E(t), n^\pm_E(t)\right)= \left( - \dot{f}(t)\pm{f^2(t)n(t)}/{2},\pm \dot{f}(t)-{f^2(t)n(t)}/{2}+n(t)\right).
$$ Consider the Riccati equation 
\begin{equation}\label{ODE1}
\begin{aligned}\frac{{\rm d} {S}_E^\pm(t)}{{\rm d}t}&= \frac{m_E^\pm(t)\pm n_E^\pm(t)}{2}(S_E^\pm)^2(t)-m_E^\pm(t)= \pm \frac{n(t)({S}_E^\pm)^2(t)}{2}  +\dot{f}(t)\mp\frac{f^2(t)n(t)}{2}.
\end{aligned}\end{equation}
Since $S_E^\pm(t)=f(t)$ is a particular  solution of   equation \eqref{ODE1}, the theory of Riccati equations yields the general solutions $$S_E^\pm(t)=f(t) \text{ or }
S_E^\pm(t)=f(t)-\frac{2e^{\mp\int m(t)\mathrm{d}t}}{\pm\int n(t)e^{\mp\int m(t)\mathrm{d}t}\mathrm{d}t+c}
$$ for all $t\in I$, where $c$ is a constant.
Consequently, 
$$\begin{aligned}
Inv^+\left(Ev^+(\bgamma),S^+_E\right)(t)=&~Ev^+(\bgamma)(t)+S_E^+(t)\bmu_E^+(t)+\frac{(S_E^+)^2(t)}{2}(Ev^+(\bgamma)(t)+\bnu_E^+(t))\\
=&~\bgamma(t)+f(t)\bnu(t)+\frac{f^2(t)}{2}(\bgamma(t)+\bmu(t))\\
&+S_E^+(t)\left(-f(t)\bmu(t)-\bnu(t)-f(t)\bgamma(t)\right) +\frac{(S_E^+)^2(t)}{2}\left(\bgamma(t)+\bmu(t)\right)\\ =&~\bgamma(t)+\left(f(t)-S^+_E(t)\right)\bnu(t)+\frac{\left(f(t)-S^+_E(t)\right)^2}{2}\left(\bgamma(t)+\bmu(t)\right).
\end{aligned}$$
Similarly, it can be  proved $$ 
Inv^-\left(Ev^-(\bgamma),S^-_E\right)(t) = \bgamma(t)+\left(f(t)-S^-_E(t)\right)\bnu(t)+\frac{\left(f(t)-S^-_E(t)\right)^2}{2}\left(\bgamma(t)-\bmu(t)\right).
$$
For the spacelike Legendre curve $\bigl({Ev}^\pm(\bgamma), -\bnu^\pm_E\bigr)$, we set $\overline{S}_E^\pm = -S_E^\mp$. Then applying Proposition \ref{prop4.10} together with the above results gives
$$
\begin{aligned}
\overline{Inv}^+\left(Ev^-(\bgamma),\overline{S}^+_E\right)(t)&=	 Inv^-\left(Ev^-(\bgamma),S^-_E\right)(t)\\
&=\bgamma(t)+\left(f(t)-S^-_E(t)\right)\bnu(t)+\frac{\left(f(t)-S^-_E(t)\right)^2}{2}\left(\bgamma(t)-\bmu(t)\right), \\
\overline{Inv}^-\left(Ev^+(\bgamma),\overline{S}^-_E\right)(t)&=	 Inv^+\left(Ev^+(\bgamma),S^+_E\right)(t)\\
&=\bgamma(t)+\left(f(t)-S^+_E(t)\right)\bnu(t)+\frac{\left(f(t)-S^+_E(t)\right)^2}{2}\left(\bgamma(t)+\bmu(t)\right).
\end{aligned}
$$
Since 
$$
f(t)-S_E^\pm(t)=0 \text{ or } f(t)-S_E^\pm(t)=\frac{2e^{\mp\int m(t)\mathrm{d}t}}{\pm\int n(t)e^{\mp\int m(t)\mathrm{d}t}\mathrm{d}t+c}
$$ 
holds for all $t\in I$, Remark \ref{re3.12} yields the conclusion that $Inv^+ (Ev^+(\bgamma),S^+_E)$, $Inv^- (Ev^-(\bgamma),S^-_E)$, $\overline{Inv}^+ (Ev^-(\bgamma),\overline{S}^+_E)$ and 
$\overline{Inv}^- (Ev^+(\bgamma),\overline{S}^-_E)$ are horocyclic parallels of $\bgamma$.	
\enD

\section{Examples}\label{S5}

\begin{example}\label{example1}{\rm
Let $(\bgamma,\bnu):(-\pi,\pi]\rightarrow\Delta_1$   be a spacelike Legendre curve given by
$$
\begin{aligned}
\bgamma(t) =\left(1+\frac{\sin^2t}{2},\frac{\sin^2t}{2},-\sin t\right), \ 
\bnu(t) =\left(\frac{\sin^2t}{2},\frac{\sin^2t}{2}-1,-\sin t\right).
\end{aligned}
$$ 
Then 
$ 
\bmu(t)=\bgamma(t)\wedge\bnu(t)=\left(\sin t,\sin t,-1\right) 
$ and
the curvature of $(\bgamma,\bnu) $ is $(m(t),n(t))=(\cos t, \cos t)$. 
Therefore, there exists a unique function $f(t)=-1$ such that $f(t)n(t)+m(t)=0$ for all $t \in (-\pi,\pi]$ and the solutions of 
$$
\frac{{\rm d}\lambda _\pm}{{\rm d}t}(t)=\mp \left(\frac{ \lambda _\pm^2(t) \cos t}{2}+\lambda _\pm(t)\cos t\right)
$$ 
are $\lambda _\pm(t)=-2 \text{ or }\lambda_\pm(t)={2c}/({e^{\pm\sin t}-c})$ for all $t\in  (\pi,\pi]$, where $c$ is a constant. 
The horocyclic evolute of $\bgamma$ is
\begin{equation}\notag
\begin{aligned}
Ev^\pm(\bgamma)(t)&=\bgamma(t)+f(t)\bnu(t)+\frac{f^2(t)}{2}(\bgamma(t)\pm\bmu(t))\\
&=\left(\frac{3}{2}\pm\frac{\sin t}{2}+\frac{\sin^2t}{4},1\pm\frac{\sin t}{2}+\frac{\sin^2t}{4},-\frac{\sin t}{2}\mp\frac{1}{2}\right).
\end{aligned} 
\end{equation}  
The horocyclic parallels of $\bgamma$ are $$P^\pm(\bgamma,-2)(t)=\left(3 \pm2\sin t  +\frac{\sin^2t}{2},2\pm2\sin t+\frac{\sin^2t}{2},- \sin t  \mp2\right)$$ and
\begin{equation}\notag
\begin{aligned}
P^\pm(\bgamma,\lambda_\pm)(t)&=\bgamma(t)+\lambda _\pm(t)\bnu(t)+\frac{ \lambda _\pm^2(t)}{2}(\bgamma(t)\pm\bmu(t))\\
&=\bigg(1+\frac{\sin^2 t}{2}+\frac{\sin^2 t}{e^{\pm\sin t}-c}+\frac{c^2 }{(e^{\pm\sin t}-c)^2} \left(2\pm2\sin t+\sin^2 t\right),\\&~~~~~~~\frac{\sin^2 t}{2}+\frac{ 2 c}{e^{\pm\sin t}-c}\left(\frac{\sin^2 t}{2}-1\right)+\frac{c^2 }{(e^{\pm\sin t}-c)^2}\left(\pm2\sin t+\sin^2 t\right),\\
&~~~~~~~-\sin t-\frac{ 2 c\sin t}{e^{\pm\sin t}-c}-\frac{c^2 }{(e^{\pm\sin t}-c)^2}\left(2\sin t \pm2\right)
\bigg). 
\end{aligned} 
\end{equation} 
\par
\begin{figure} 
\centering
\includegraphics[width = 9.5cm]{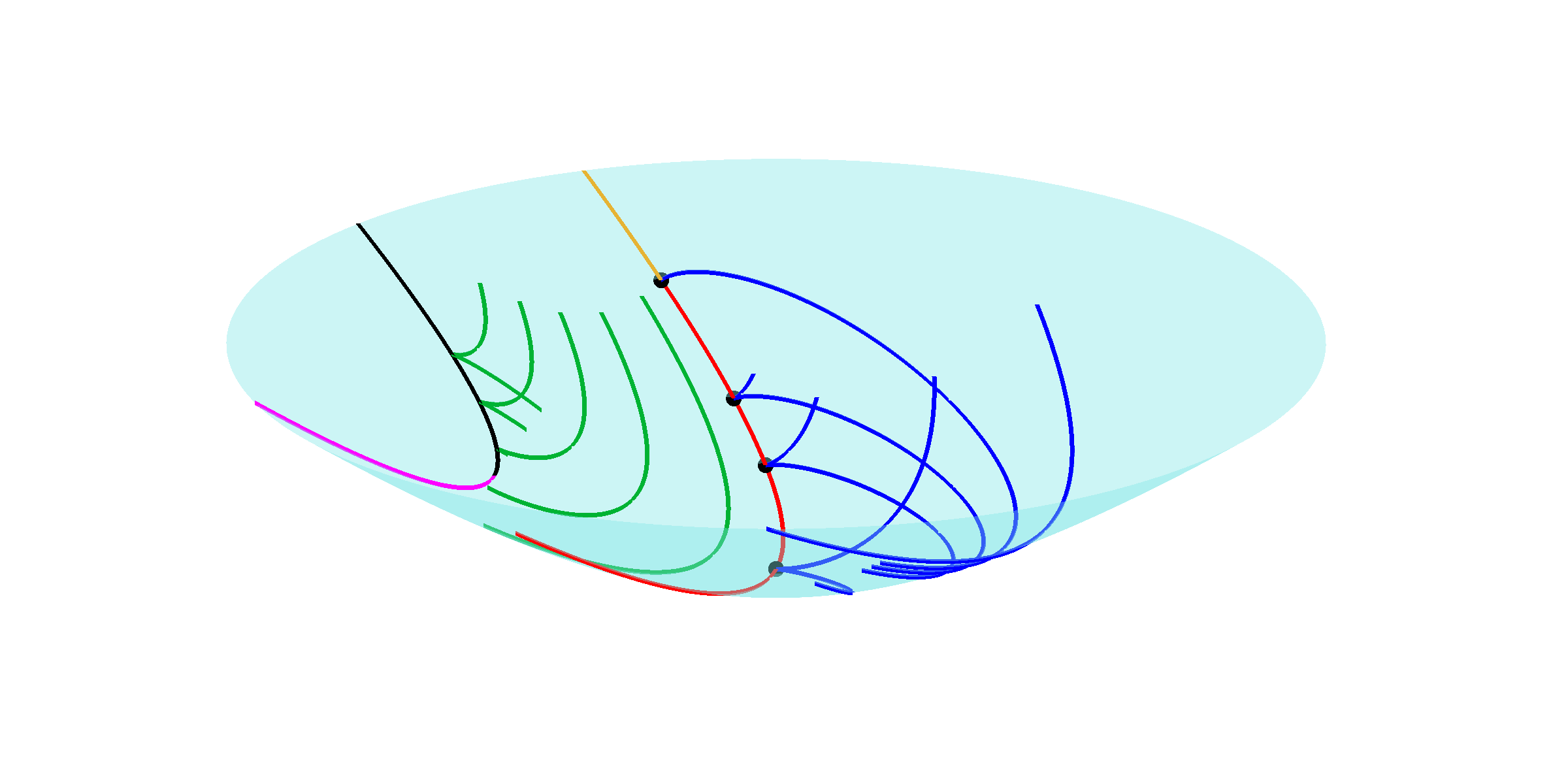}
\caption{A spacelike frontal $\bgamma$ (red curve) and its horocyclic involutes   with $\overline{c}=-1,-e^{-1},-e^{-\sqrt{2}},-e^{-2}$, $\widetilde{c}=0,1/2,\sqrt{2}/2,1$ and $s_+(t)=-1$ (blue curves and black points), horocyclic evolutes (magenta and black curves),   horocyclic parallels (green and yellow curves) with $c=0.4,0.6,0.8,1$ and $\lambda_\pm(t)=-2$ in Example \ref{example1}. }
\label{figure1}
\end{figure}

Furthermore, we consider the ordinary differential equations
\begin{equation} \notag
\frac{{\rm d}s_+}{{\rm d}t}(t)=  (\cos t) s_+^2(t)-\cos t \text{ and }   
\frac{{\rm d}s_-}{{\rm d}t}(t)= -\cos t.
\end{equation}
The solutions of the equations are given respectively by 
$$ s_+(t)= -1
 \text{ or } s_+(t)=\frac{1+\overline{c}e^{2\sin t }}{1-\overline{c}e^{2\sin t }},$$ and $s_-(t)=-\sin t+\widetilde{c}$ for all $t \in  (-\pi,\pi]$, where $\overline{c}$ and $\widetilde{c}$ are  constants.
Hence,  the horocyclic involutes of $\bgamma$ are 
$$
\begin{aligned}
Inv^+(\bgamma,-1)(t)&=\left(\frac{3}{2} -\sin t+\sin^2t,-\frac{1}{2}-\sin t+\sin^2t,1-2\sin t\right),\\
Inv^+(\bgamma,s_+)(t)&=\bgamma(t)+s_+(t)\bmu(t)+\frac{s_+^2(t)}{2}(\bgamma(t)+\bnu(t))\\
&=\bigg(1+\frac{\sin^2t}{2}+\frac{\left(1+\overline{c}e^{2\sin t}\right)\sin t}{1-\overline{c}e^{2\sin t}}+\frac{\left(1+\overline{c}e^{2\sin t}\right)^2(1+\sin^2t)}{2\left(1-\overline{c}e^{2\sin t}\right)^2},\\
&~~~~~~\frac{\sin^2t}{2}+\frac{\left(1+\overline{c}e^{2\sin t}\right)\sin t}{1-\overline{c}e^{2\sin t}}+\frac{\left(1+\overline{c}e^{2\sin t}\right)^2(\sin^2t-1)}{2\left(1-\overline{c}e^{2\sin t}\right)^2},\\
&~~~~~~-\sin t-\frac{ 1+\overline{c}e^{2\sin t} }{1-\overline{c}e^{2\sin t}}-\frac{\left(1+\overline{c}e^{2\sin t}\right)^2 \sin t}{ \left(1-\overline{c}e^{2\sin t}\right)^2}\bigg)
\end{aligned}
$$ 
and
$$
\begin{aligned}
Inv^-(\bgamma,s_-)(t)&=\bgamma(t)+s_-(t)\bmu(t)+\frac{s_-^2(t)}{2}(\bgamma(t)-\bnu(t)) =\left(1+\frac{\widetilde{c}^2 }{2},\frac{\widetilde{c}^2 }{2},-\widetilde{c}\right).
\end{aligned}
$$ 
 
If $\overline{c}=-1$, since $s_+(0)=0$, $m(0)=1$ and $m(0)+n(0)=2$, we can get $ Inv^+(\bgamma,s_+) $ has a $(2,3)$-cusp at $t=0$  by Corollary \ref{cor4.8}, see Figure \ref{figure1}.
}
\end{example}

\begin{example}\label{example2}{\rm 

Let $(\bgamma,\bnu):(-\pi,\pi]\rightarrow\Delta_1$ be a spacelike Legendre curve given by
$$
\begin{aligned}
\bgamma(t)&= \left(1+\frac{g^2(t)}{2}\right)\bx_0(t)+g(t)\bx_1(t)+\frac{g^2(t)}{2}\bx_2(t),\\
\bnu(t)&=\left(1-\frac{g^2(t)}{2}\right)\bx_2(t)-g(t)\bx_1(t)-\frac{g^2(t)}{2}\bx_0(t),
\end{aligned}
$$
where 
\begin{equation}\notag
\begin{aligned}
\bx_0(t)&=	\left(\sqrt{1+\cos^6t+\sin^6t},\cos^3t,\sin^3t\right),\\	
\bx_1(t)&=\frac{\left(\sin t\cos t\sqrt{1+\cos^6t+\sin^6t},\sin t(1+\cos^4t),\cos t(1+\sin^4t)\right)}{\sqrt{1+\cos^2t\sin^2t}},\\
\bx_2(t)&=\frac{\left( \sin^2t-\cos^2t,-\cos t\sqrt{1+\cos^6t+\sin^6t},\sin t\sqrt{1+\cos^6t+\sin^6t}\right)}{\sqrt{1+\cos^2t\sin^2t}},\\
g(t)&=-\frac{3\sin t\cos t(1+\sin^2t\cos^2t)^{3/2}}{3\sin^4t\cos^4t+3\sin^2t\cos^2t-1-\sin^6t-\cos^6t}.
\end{aligned}		
\end{equation}
 
\begin{figure}[t] 
	\begin{minipage}{0.45  \textwidth}
		\centering
		\includegraphics[width = 9cm]{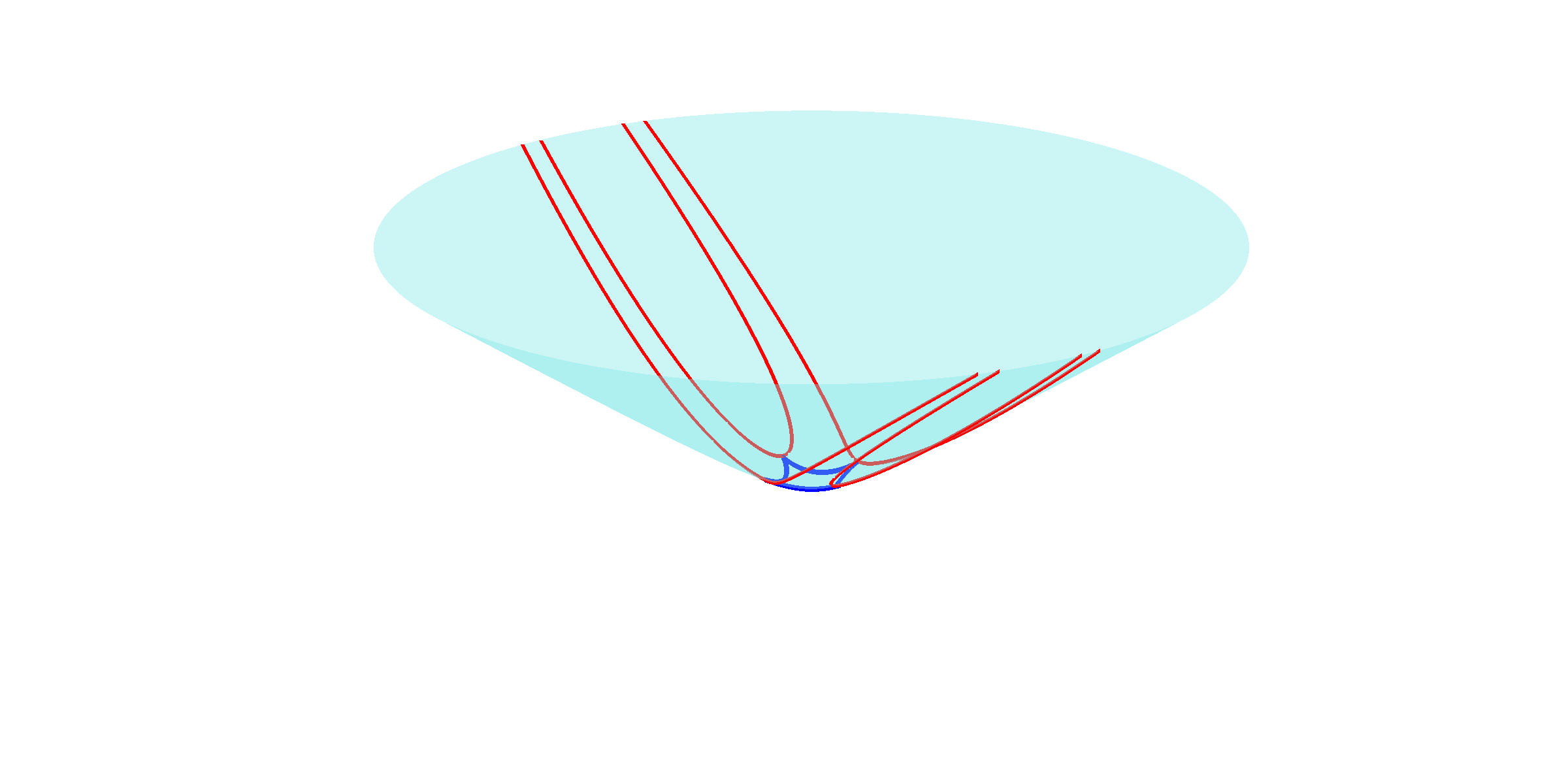}
		\caption{A spacelike frontal $\bgamma$ (red curve) and its horocyclic involute    (blue curve) in Example  \ref{example2}.}
		\label{figure2}
	\end{minipage}
	\qquad
	\begin{minipage}{0.45 \textwidth}
		\centering
		\includegraphics[width = 9cm]{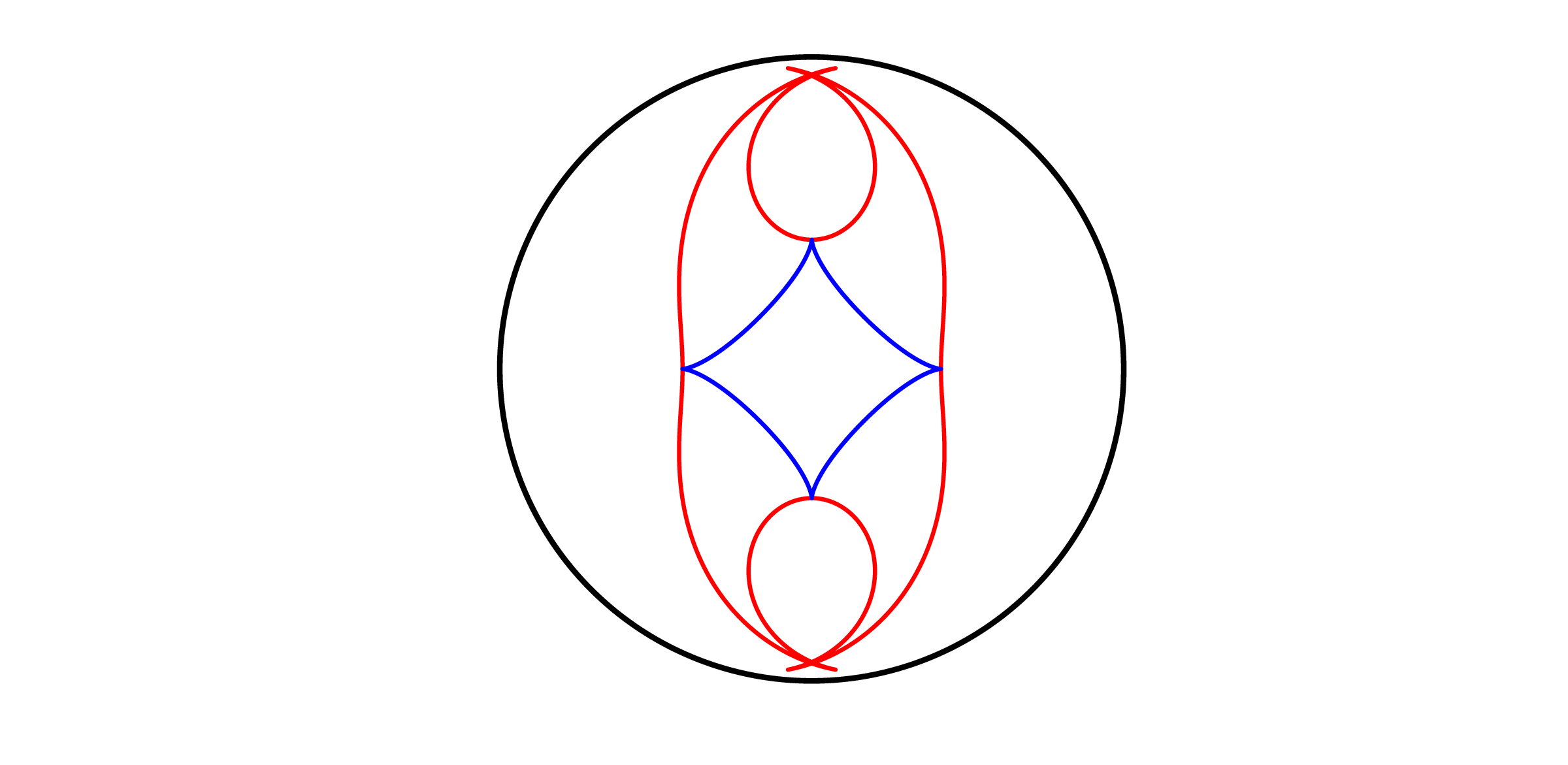}
		\caption{A spacelike frontal $\bgamma$ (red curve) and its  horocyclic involute (blue curve) projected to Poincar\'{e} 2-disc in Example \ref{example2}.}
		\label{figure3}
	\end{minipage}
\end{figure}		
\begin{figure} 
	\begin{minipage}[t]{0.45  \textwidth}
		\centering
		\includegraphics[width = 9cm]{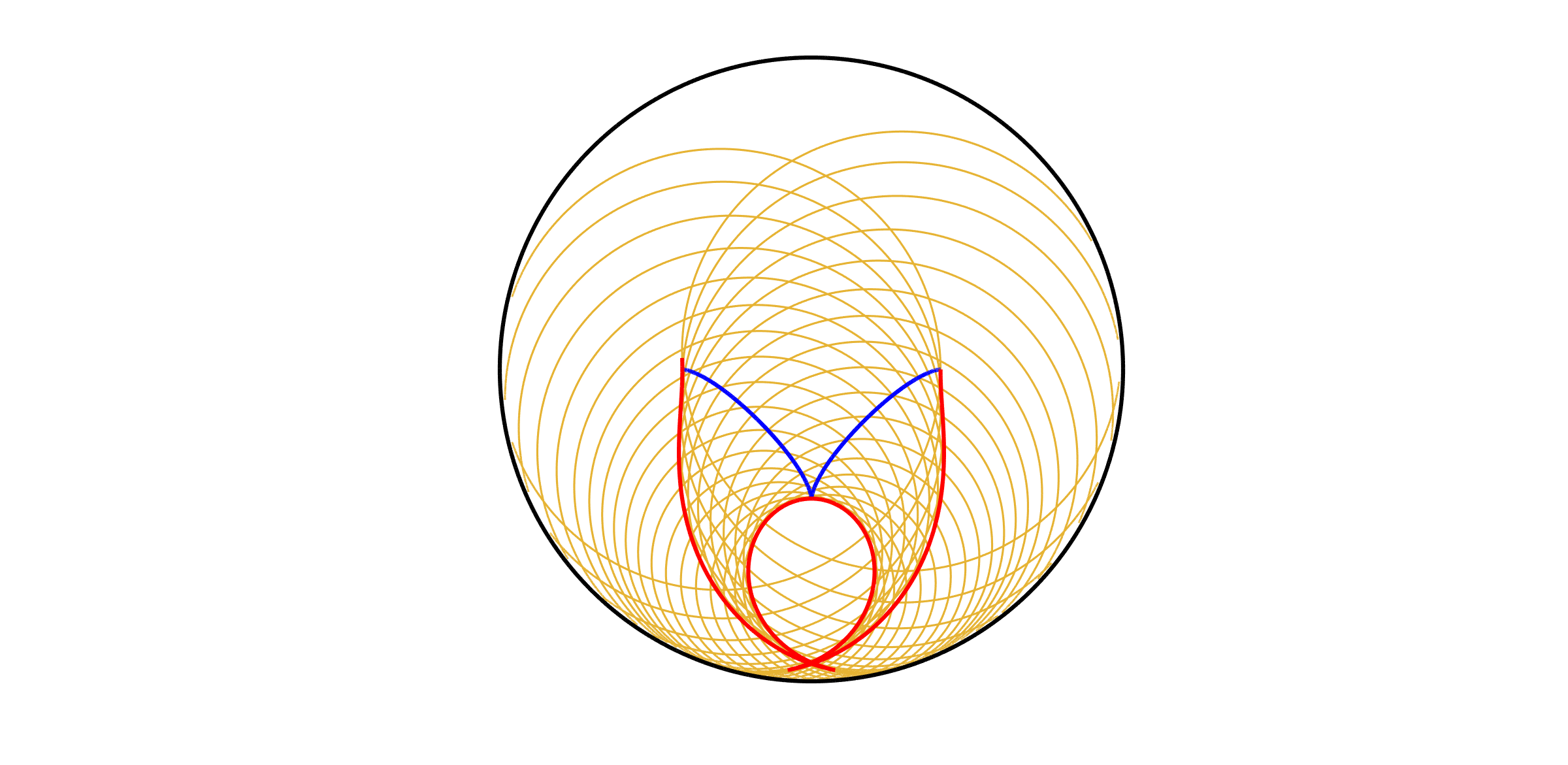}
		\caption{A horocyclic involute (blue curve) of $\bgamma$ as a normal envelope of the tangent horocycles of $\bgamma$ (red curve) in Example  \ref{example2}.}
		\label{figure4}
	\end{minipage}
	\qquad
	\begin{minipage}[t]{0.45  \textwidth}
		\centering
		\includegraphics[width = 9cm]{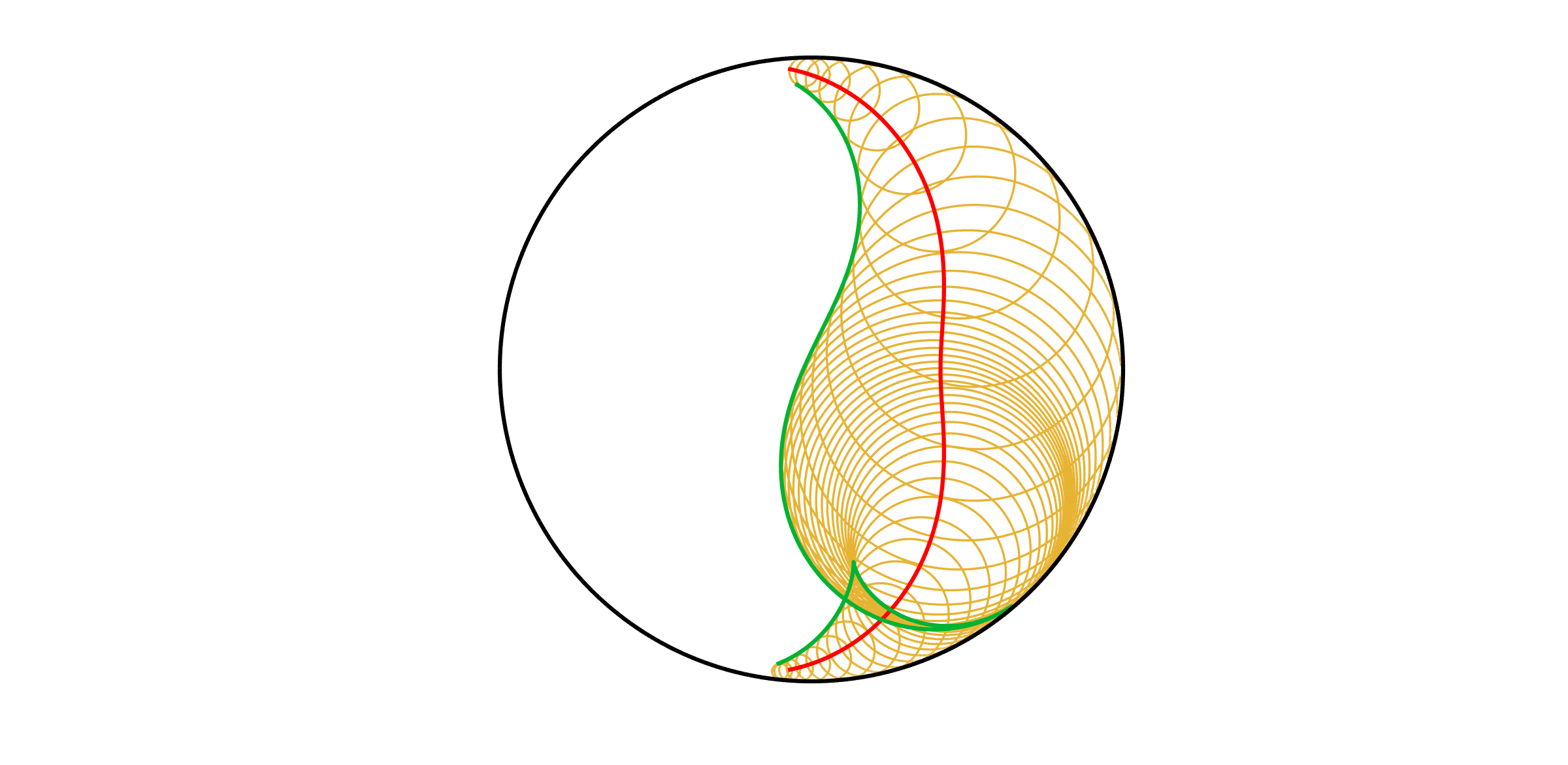}
		\caption{A horocyclic evolute (green curve) of $\bgamma$ as an  envelope of the normal horocycles of $\bgamma$ (red curve) in Example \ref{example2}.}
		\label{figure5}
	\end{minipage}
\end{figure}
Then 
$ 
\bmu(t)=\bgamma(t)\wedge\bnu(t)= -\bx_1(t)-g(t)(\bx_0(t)+\bx_2(t)).
$ 
The curvature of $(\bgamma,\bnu) $ is given by 
$$
(m(t),n(t))=\left( -\dot{g}(t)+\frac{g^2(t)e(t)}{2},\dot{g}(t)-\frac{g^2(t)e(t)}{2}+e(t)\right),
$$
where 
$$
e(t)=\frac{3\sin^4t\cos^4t+3\sin^2t\cos^2t-1-\sin^6t-\cos^6t}{(1+\sin^2t\cos^2t)\sqrt{1+\cos^6t+\sin^6t}}.
$$
The function $n(t)$ has two approximate zeros at $-t_0$ and $t_0$,  where $t_0=0.7986$ and the limit of $ n $ fails to exist at these points. However,
there exists a smooth function $f: (-t_0,t_0)\rightarrow \R$ such that $f(t)=-m(t)/n(t)$ for all $t\in (-t_0,t_0)$.     Therefore,   the horocyclic evolute of $\bgamma$ in $(-t_0,t_0)$ is 
$$
Ev^\pm(\bgamma)(t)=\bgamma(t)-\frac{m(t)}{n(t)}\bnu(t)+\frac{m^2(t)}{2n^2(t)}(\bgamma(t)\pm\bmu(t)),
$$
see Figure \ref{figure5}.
Consider the Riccati equation
\begin{equation}\label{ODE2}
\frac{{\rm d}s_+}{{\rm d}t}(t)= \frac{m(t)+ n(t)}{2}s_+^2(t)-m(t).
\end{equation}
It can be verified that $ s_+(t)=g(t) $ is a solution of equation \eqref{ODE2}.
Hence, the corresponding horocyclic involute of $\bgamma$ is given by
$$
\begin{aligned}
Inv^+(\bgamma,s_+)(t)&=\bgamma(t)+g(t)\bmu(t)+\frac{g^2(t)}{2}(\bgamma(t)+\bnu(t))  =\left(\sqrt{1+\cos^6t+\sin^6t},\cos^3t,\sin^3t\right). 
\end{aligned}
$$
If $t_1=-\pi/2,0,\pi/2,\pi$,  the fact $s_+(t_1)=g(t_1)=0$ implies that  $t_1$ is a singular point of $ Inv^+(\bgamma,s_+)$ by Corollary \ref{cor4.8} (1). 
Take $t_1=0$ as an example, since $m(0)+n(0)=-\sqrt{2}\neq0$, $m(0)=-\dot{g}(0)=-3/2\neq0$, we obtain that $ Inv^+(\bgamma,s_+)$ has a (2,3)-cusp at $t_1=0$ by Corollary \ref{cor4.8}. 
A similar discussion leads to the conclusion that $Inv^+(\bgamma,s_+)$ is also diffeomorphic to $(2,3)$-cusp at $t_1=-\pi/2,\pi/2,\pi$, see Figures \ref{figure2}, \ref{figure3}, \ref{figure4}.
}
\end{example}


Nozomi Nakatsuyama, 
\\
Muroran Institute of Technology, Muroran 050-8585, Japan,
\\
E-mail address: 25096009b@muroran-it.ac.jp
\\
\\
Masatomo Takahashi, 
\\
Muroran Institute of Technology, Muroran 050-8585, Japan,
\\
E-mail address: masatomo@muroran-it.ac.jp
\\
\\
Anjie Zhou,
\\
School of Mathematics and Statistics, Northeast Normal University, Changchun 130024, P. R. China
\\
E-mail address: zhouaj882@nenu.edu.cn

\end{document}